\documentclass[a4paper, 11pt]{scrartcl}
\pdfoutput=1

\usepackage[utf8]{luainputenc}
\usepackage[english]{babel}
\usepackage{csquotes}

\usepackage[a4paper, top=0.51in, bottom=0.79in, left=0.999in, right=0.99in]{geometry}

\usepackage[plainpages=false,pdfpagelabels,hidelinks]{hyperref}

\usepackage[%
  backend=bibtex,bibencoding=ascii,
  style=authoryear-comp, dashed=false,
  firstinits=true, uniquename=init, 
  natbib=true,
  url=true,
  doi=true,
  isbn=false,
  backref=false,
  maxnames=99,
  ]{biblatex}
\usepackage{amsmath}
\usepackage{amssymb}
\usepackage{commath}
\usepackage{mathtools}

\usepackage{amsthm}
\theoremstyle{plain}
  \newtheorem{thm}{Theorem}
  \newtheorem{lem}[thm]{Lemma}

\usepackage{graphicx}
\usepackage[small]{caption}
\usepackage{subcaption}

\usepackage{booktabs}
\usepackage{rotating}

\usepackage{xparse}

\renewcommand{\vec}[1]{\underline{#1}}
\NewDocumentCommand{\mat}{mo}{%
  \IfValueTF{#2}{%
    \underline{\underline{#1}}{#2}
  }{%
    \underline{\underline{#1}}\,
  }%
}
\newcommand{\vecfnum}{\vec{f}^\mathrm{num}}
\newcommand{\fnum}{f^\mathrm{num}}

\DeclarePairedDelimiter{\diagfences}{(}{)}
\newcommand{\diag}{\operatorname{diag}\diagfences}

\newcommand{\scp}[2]{\left\langle{#1, #2}\right\rangle}

\renewcommand{\epsilon}{\varepsilon}
\renewcommand{\phi}{\varphi}
\newcommand{\N}{\mathbb{N}}
\newcommand{\R}{\mathbb{R}}

\usepackage{ifthen}
\newboolean{numeric_citation}
\setboolean{numeric_citation}{false}

\ifthenelse{\boolean{numeric_citation}}{
  \newcommand{\citex}[2]{#1}
}{
  \newcommand{\citex}[2]{#2}
}

\title{Enhancing stability of correction procedure via reconstruction using
       summation-by-parts operators II: Modal filtering}
\author{Jan Glaubitz, Hendrik Ranocha, Philipp Öffner, Thomas Sonar}
\titlehead{\includegraphics[width=0.5\textwidth]{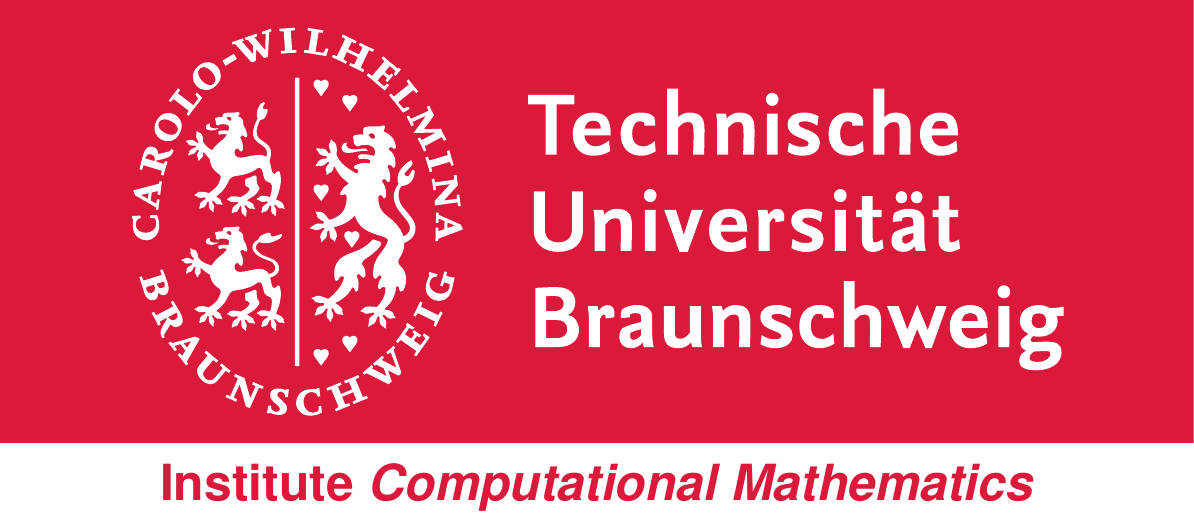}}
\date{June 2, 2016}

\hypersetup{pdfauthor={Hendrik Ranocha}}
\hypersetup{pdftitle={Modal filtering for CPR using SBP operators}}
\hypersetup{pdfkeywords={hyperbolic conservation laws, high order methods,
summation-by-parts, SBP, flux reconstruction, FR, lifting collocation penalty,
LCP, correction procedure via reconstruction, CPR, entropy stability,
modal filtering}}

\begin{document}
  \maketitle
  
  \begin{abstract}
  A recently introduced framework of semidiscretisations for hyperbolic
  conservation laws known as correction procedure via reconstruction (CPR,
  also known as flux reconstruction) is considered in the extended setting
  of summation-by-parts (SBP) operators using simultaneous approximation terms
  (SATs). This reformulation can yield stable semidiscretisations for linear
  advection and Burgers' equation as model problems. In order to enhance these
  properties, modal filters are introduced to this framework.
  As a second part of a series, the results of Ranocha, Glaubitz, Öffner, and
  Sonar (\emph{Enhancing stability of correction procedure via reconstruction using
  summation-by-parts operators I: Artificial dissipation}, 2016)
  concerning artificial dissipation / spectral viscosity are extended,
  yielding fully discrete stable schemes. Additionally, a new adaptive strategy
  to compute the filter strength is introduced and different possible applications
  of modal filters are compared both theoretically and numerically.
\end{abstract}

  \section{Introduction}

Numerous phenomena in the field of engineering as well as natural sciences are
modelled by partial differential equations (PDEs) and in particular by 
hyperbolic conservation laws.
While being of vital significance, especially hyperbolic conservation laws are
highly delicate from a mathematical point of view and often exact solutions are
not known. 
Therefore, the study of numerical solutions is an important and interesting field
of research.

%
In the last decades, great effort was made to develop special numerical schemes
for these PDEs. 
Often the method of lines using a simple \emph{Runge-Kutta} (RK) scheme in time 
is adapted \citep{gottlieb1998total,ketcheson2008highly}, 
where far most attention is given
to suitable space discretisations. 
These can be obtained by well-known low order methods  as well as high order 
ones, e.g. \emph{discontinuous Galerkin} (DG), \emph{spectral difference} (SD), 
recent \emph{flux reconstruction} (FR) or  
\emph{correction procedure via reconstruction} (CPR) methods. 
While the latter ones can provide highly desirable accuracy for smooth solutions, 
they are originally far less robust in the case of discontinuities than their low 
order counterparts.
However, using \emph{summation-by-parts} (SBP)
operators along with \emph{simultaneous-approximation-terms} (SATs) provably 
stable and conservative space discretesations in the FD framework were obtained. 
See for instance the review articles \citex{}{by}
\cite{svard2014review, nordstrom2015new, fernandez2014review}.
In a similar manner their ideas have been extended to stable space 
discretisations for the CPR method \citex{in}{by} \citet{ranocha2016summation}.
The FR method along with the CPR method as its generalisation provide a 
cell-wise spectral approach on a decomposition of the underlying domain, 
using correction terms at the cell boundaries quite similar to the SATs from 
the FD framework.
Hence, by adapting analytical counterparts of the SBP operators using polynomial 
bases, also stable space discretisations for the FR and CPR method could be 
obtained.

Unfortunately, these stable space semidiscretisations alone cannot provide 
stable schemes in the fully discrete setting. 
Coming back to the decoupled time integration, 
already the explicit Euler method 
yealds an additional error after every time step. 

%
This error might become most destructive when schemes using polynomial 
approximations like the CPR method are applied.
Here, the Gibbs phenomenon will arise for solutions with jump 
discontinuities, and the method becomes unstable.
There are several well-known solutions to enhance stability
of a numerical scheme. 
In this work only the artificial dissipation or spectral viscosity method
and modal filtering by certain exponential filters will be addressed. 
Artificial dissipation / spectral viscosity can be interpreted as 
discrete counterparts of the vanishing viscosity method.
The latter one was originally used to prove the existence of
entropy solutions by adding an iteratively decreasing viscosity term to the 
right hand
side of the conservation law and constructing the $L_1$ limit of solutions to
the resulting parabolic equation \citep{lax1973hyperbolic}. 
Since then, appropriate discretisations of such viscosity terms can be used to
stabilise numerical schemes for hyperbolic conservation laws.
If this discretisation is done in the FD framework, 
the corresponding method is often labelled 
\emph{artificial dissipation} (AD) 
\citep{mattsson2004stable, nordstrom2006conservative}
Using a method which approximates the solution polynomials however, 
the resulting method is then usually named  
\emph{spectral viscosity} (SV) 
\citep{ma1998chebyshev1, ma1998chebyshev2, maday1989analysis}.

On the other hand, the idea of modal filtering exclusively adapts to 
such methods using series expansions \citep{hesthaven2008filtering}.
This has already be done for DG methods \citex{in}{by} 
\citet{meister2012application, meister2013extended} 
and for the SD method \citex{in}{by} \citet{glaubitz2016application}.
%
At least for certain exponential filters, there is an equivalence to the
spectral viscosity method, if a splitting in the time integration is used.
In this work, both approaches are applied to the new CPR method 
using SBP operators.
While this has already been done in various well-known schemes, 
this work is the first to do so in this particular recent one. 
Furthermore, the framework of SBP operators enables the study of 
different stabilisation approaches to the
numerical solution and in particular to already mentioned errors from 
time integration schemes like the explicit Euler method.
The authors make use of this by deriving completely new adaptive filtering 
techniques, which have not been published before to their best knowledge.
While AV is addressed in the first part \citep{ranocha2016enhancing} 
of this series, this one focuses on modal filtering.
Doing so, the rest is organised as follows:

In section 2 a brief introduction to the correction procedure via 
reconstruction using summation-by-parts operators is given. 
This recent method is proven to be conservative and stable for the problem of
linear advection as well as for the non-linear Burgers' equation
using a skew-symmetric formulations for the discretisation.
One should note however that the latter modal filtering approaches along with
their analysis just require a stable space-discretisation using 
a polynomial approximation of the solution and thus are completely 
independent from this specific one.

This is stressed in the third section, when the idea of spectral viscosity
and its connection to modal filtering is addressed.
While part I of this work focused mainly on artificial dissipation, 
this part focuses on the idea of modal filtering and its strong 
relation by a certain operator splitting in time.

Naturally the question arises, when exactly modal filtering 
should be applied.
Section 4 discusses the straight forward option of filtering in the sense
of the above mentioned operator splitting, the fifth and sixth section
address the potential of filtering the time derivative 
and filtering the solution.

The work is then completed by the authors' conclusions and 
quite a few remarks about further research in the seventh section.

  \section{Correction procedure via reconstruction using summation-by-parts operators}

In this section, a brief description of CPR methods using SBP operators will be
given. Therefore, a one-dimensional scalar conservation law
\begin{equation}
\label{eq:scalar-CL}
  \partial_t u + \partial_x f(u) = 0,
\end{equation}
equipped with an adequate initial condition will be used. For simplicity,
periodic boundary conditions (or a compactly supported initial condition) will
be assumed.

The domain is partitioned into subdomains, which are mapped diffeomorphically
onto a reference element. Since only one spatial dimension is considered, these
intervals can be mapped by an affine linear transformation. The following
computations are done in this standard domain.

In the setting described \citex{in}{by} \citet{ranocha2015extended}, the solution is represented
as an element of a finite dimensional Hilbert space of functions on the volume.
With respect to a chosen basis, the scalar product approximating the $L_2$ scalar
product is represented by a matrix $\mat{M}$ and the derivative (divergence) by
$\mat{D}$. Additionally, functions on the boundary (consisting of two points in
this one dimensional case) are elements of another finite dimensional Hilbert
space with appropriate basis. The restriction of functions on the volume to the
boundary is represented by a (rectangular) matrix $\mat{R}$ and integration
with respect to the outer normal by $\mat{B}$. Finally, the operators have to
fulfil the SBP property
\begin{equation}
\label{eq:SBP}
  \mat{M} \mat{D} + \mat{D}[^T] \mat{M}
  = \mat{R}[^T] \mat{B} \mat{R}
\end{equation}
as compatibility condition in order to mimic integration by parts
\begin{equation}
  \vec{u}^T \mat{M} \mat{D} \vec{v} + \vec{u}^T \mat{D}[^T] \mat{M} \vec{v}
  \approx
  \int_\Omega u \, \partial_x v + \int_\Omega \partial_x u \, v
  = u \, v \big|_{\partial \Omega}
  \approx
  \vec{u}^T \mat{R}[^T] \mat{B} \mat{R} \vec{v}.
\end{equation}
Additional ingredients of CPR methods are numerical fluxes (Riemann solvers) $\fnum$,
computing a common flux on the boundary using values from both neighbouring
elements, and a correction step, which \citex{Ranocha, Öffner, and Sonar}{} \citet{ranocha2016summation} formulated
as an SAT, similarly to the weak enforcement of boundary condition in finite
difference SBP schemes.

In the following, polynomial bases will be used for functions on the volume,
either nodal Gauß-Legendre / Lobatto-Legendre or modal Legendre bases.
Multiplication of functions on the volume will be conducted pointwise for
nodal bases or exactly, followed by an $L_2$ projection, for modal bases. The
resulting multiplication operators are written with two underlines, e.g.
$\mat{u}$ represents multiplication with the polynomial given by $\vec{u}$.

\subsection{Linear advection}

A simple example for a scalar conservation law \eqref{eq:scalar-CL} is given by
the linear advection equation with constant velocity
\begin{equation}
\label{eq:lin-adv}
  \partial_t u + \partial_x u = 0.
\end{equation}
The semidiscretisation using the canonical of the correction procedure can be
written as
\begin{equation}
\label{eq:lin-adv-semidisc}
  \partial_t \vec{u}
  = - \mat{D} \vec{u}
    - \mat{M}[^{-1}] \mat{R}[^T] \mat{B} \left(
        \vecfnum - \mat{R} \vec{u}
      \right)
\end{equation}
in the reference domain. It is conservative across elements and stable if an
appropriate numerical flux is applied, see inter alia Theorem 5 of
\citet{ranocha2016summation}.

\subsection{Burgers' equation}

As a nonlinear model, Burgers' equation
\begin{equation}
\label{eq:Burgers}
  \partial_t u + \partial_x \frac{u^2}{2} = 0
\end{equation}
is more difficult to handle correctly, since discontinuities may develop in
finite time. However, a conservative (across elements) and stable (with respect
to the discrete norm induced by $\mat{M}$) semidiscretisation can be obtained by
the application of a skew-symmetric form
\begin{equation}
\label{eq:Burgers-semidisc}
  \partial_t \vec{u}
  =
  - \frac{1}{3} \mat{D} \mat{u} \vec{u}
  - \frac{1}{3} \mat{u}{^*} \mat{D} \vec{u}
  + \mat{M}[^{-1}] \mat{R}[^T] \mat{B} \left(
      \vecfnum
      - \frac{1}{3} \mat{R} \mat{u} \vec{u}
      - \frac{1}{6} \left( \mat{R} \vec{u} \right)^2
    \right),
\end{equation}
if the numerical flux is adequate, see inter alia Theorem 2 of
\citet{ranocha2015extended}.

  \section{Spectral viscosity and modal filtering}
\label{sec:SV-and-filtering}

In this section, spectral viscosity will be briefly introduced as a means to
improve stability properties not influencing conservation across elements.
Using the Sturm-Liouville operator associated with Legendre polynomials for
artificial dissipation, the continuous properties for polynomial bases are
investigated in the discrete setting and modal filters are derived by a
classical operator splitting approach.

\subsection{Spectral viscosity in the continuous and semidiscrete setting}

A spectral viscosity extension of a scalar conservation law \eqref{eq:scalar-CL}
can be written as
\begin{equation}
\label{eq:scalar-CL-RHS}
  \partial_t u(t,x) + \partial_x f(u(t,x))
  = (-1)^{s+1} \epsilon \left( \partial_x a(x) \partial_x \right)^{s} u(t,x),
\end{equation}
introducing the \emph{strength} $\epsilon \geq 0$, the order $s \in \N$, and a
suitable polynomial $a \colon \R \to \R$, fulfilling $a \big|\Omega \geq 0$
and $a \big|_{\partial \Omega} = 0$. Then, conservation across the element
$\Omega$ is preserved and the artificial dissipation on the right hand side 
enforces an additional decay of the $L_2$ norm of the solution $u$, as described
inter alia \citex{in}{by} \citet{ranocha2016enhancing} in section 3.1.

In order to enforce the same properties in the semidiscrete setting, the
artificial dissipation should be discretised as
\begin{equation}
\label{eq:RHS-smart}
  (-1)^{s+1} \epsilon \left( - \mat{M}[^{-1}] \mat{D}[^T] \mat{M} \mat{a} \mat{D}
  \right)^s \vec{u}
  =
  - \epsilon \left( \mat{M}[^{-1}] \mat{D}[^T] \mat{M} \mat{a} \mat{D} \right)^{s} \vec{u},
\end{equation}
see inter alia section 3.2 of \citet{ranocha2016enhancing}.
Additionally, they prove
\begin{lem}[Lemma 2 of \citet{ranocha2016enhancing}]
\label{lem:eigenvalues}
  The discrete viscosity operator $- \mat{M}[^{-1}] \mat{D}[^T] \mat{M} \mat{a} \mat{D} u$
  of the right hand side \eqref{eq:RHS-smart}
  \begin{itemize}
    \item 
    has the correct eigenvalues $- n (n+1)$ with eigenvectors given by the
    Legendre polynomials $\phi_n, n \in \set{0, \dots, p}$ if a modal Legendre
    or a nodal Gauß-Legendre basis is used. These eigenvalues and eigenvectors
    are the same as the ones of the continuous operator $\partial_x (1-x^2) \partial_x$.
    
    \item
    has the correct eigenvalues $- n (n+1)$ with eigenvectors given by the
    Legendre polynomials $\phi_n, n \in \set{0, \dots, p-1}$, if a nodal
    Lobatto-Legendre basis is used.
    The eigenvalue for $\phi_p$ is non-positive and bigger than $- p (p+1)$,
    i.e. the artificial viscosity operator yields less dissipation of the highest
    mode compared to the exact value, obtained by modal Legendre and nodal
    Gauß-Legendre bases.
  \end{itemize}
\end{lem}

\subsection{Discrete setting}
\label{sec:discrete-setting}

In order to discretise the augmented conservation law \eqref{eq:scalar-CL-RHS}
\begin{equation}
  \partial_t u + \partial_x f(u) = \epsilon \, \partial_x a(x) \partial_x u
\end{equation}
with $a(x) = 1-x^2$, discretisations for both space
and time have to be chosen. Here, a conservative and stable semidiscretisation
is obtained as a CPR method using SBP operators. Thus, the solution $u$ is a
polynomial in each element and can therefore be written as a linear combination
of Legendre polynomials, i.e. $u = \sum_{i=0}^p u_i \phi_i$. In this basis of
its eigenvectors, the right hand side
$- \epsilon \, \mat{M}[^{-1}] \mat{D}[^T] \mat{M} \mat{a} \mat{D} u$
is represented by a diagonal matrix with the eigenvalues according to Lemma
\ref{lem:eigenvalues} on the diagonal.

There are at least two choices for a time discretisation: Discretise both 
$-\partial_x f(u)$ and $\epsilon \, \partial_x a(x) \partial_x u$ together or use
some operator splitting method. The first choice if investigated \citex{in}{by}
\citet{ranocha2016enhancing}. Here, an explicit Euler method with operator
splitting of first order is applied, resulting in
\begin{equation}
  \vec{u} \mapsto \tilde{\vec{u}}_+ \mapsto \vec{u}_+ := \mat{F} \tilde{\vec{u}}_+,
\end{equation}
where $\tilde{\vec{u}}_+$ is obtained by one time step of a discretisation of
$- \partial_x f(u)$ and the linear operator $\mat{F}$ describes one time step of
the linear ordinary differential equation
\begin{equation}
  \od{}{t} \vec{u} = - \epsilon \left( 
    \mat{M}[^{-1}] \mat{D}[^T] \mat{M} \mat{a} \mat{D} \right)^s \vec{u},
\end{equation}
i.e.
\begin{equation}
  \mat{F} = \exp \left[ - \epsilon \left( 
    \mat{M}[^{-1}] \mat{D}[^T] \mat{M} \mat{a} \mat{D} \right)^s
                    \Delta t \right].
\end{equation}
Since the Legendre polynomials $\phi_n$ are eigenvectors of the viscosity operator
with eigenvalues $\lambda_n^s$, $\mat{F}$ is diagonal in this basis
\begin{equation}
  \mat{F} = \diag{ \exp \left[ -\epsilon \lambda_n^s \Delta t \right]_{n=0}^p },
\end{equation}
where the correct (continuous) eigenvalues are $\lambda_n = n (n+1)$ and the
discrete eigenvalues are given in Lemma \ref{lem:eigenvalues}.
Thus, in a modal Legendre basis, a first order operator splitting of
\eqref{eq:scalar-CL-RHS} can be written as a time step for the conservation
law with vanishing right hand side \eqref{eq:scalar-CL}, followed by an
application of the modal filter $\mat{F}$. In order to use the correct eigenvalues,
$\mat{F}$ can be obtained by a simple change of basis if a nodal basis is used.
However, this idea of modal filtering can also be applied in other ways. Some
possibilities are
\begin{enumerate}
  \item
  \label{itm:operator-splitting}
  Use the above mentioned operator splitting
  \begin{enumerate}
    \item Time step $\Delta t$ for $\partial_t u + \partial_x f(u) = 0$ yielding
          $\tilde u_+$.
    \item Filter $\tilde u_+ \mapsto u_+ := F \tilde u_+$.
  \end{enumerate}
  
  \item
  \label{itm:filter-derivative}
  Apply the filter to the time derivative:
  $ \partial_t u + F \, \partial_x f(u) = 0$.
  
  \item
  \label{itm:filter-u}
  Apply the filter to the solution used to compute the time derivative:
  $ \partial_t u + \partial_x f( F u ) = 0$.
\end{enumerate}

Possibility \ref{itm:operator-splitting} shall yield a stable method if enough
dissipation is added, since this holds for the semidiscretisation. However,
all three possibilities will be investigated in the following sections.

  \section{Operator splitting}
\label{sec:operator-splitting}

In this section, the operator splitting approach, i.e. possibility
\ref{itm:operator-splitting} of section \ref{sec:discrete-setting}, will be
considered.

\subsection{Semidiscrete and discrete estimates}

As described in the previous section, the discrete viscosity operator
$- \epsilon \left( \mat{M}[^{-1}] \mat{D}[^T] \mat{M} \mat{a} \mat{D} \right)^s$
on the right hand side yields a conservative and stable semidiscretisation of
the scalar conservation law with additional viscosity \eqref{eq:scalar-CL-RHS}
\begin{equation}
  \partial_t u + \partial_x f(u)
  = (-1)^{s+1} \epsilon \left( \partial_x a \partial_x \right)^s u.
\end{equation}
Investigating conservation across elements while applying an explicit Euler
method $\vec{u} \mapsto \vec{u}_+ := \vec{u} + \Delta t \, \partial_t \vec{u}$
as time discretisation results in
\begin{equation}
  \vec{1}^T \mat{M} \vec{u}_+
  = \vec{1}^T \mat{M} \vec{u} + \Delta t \vec{1}^T \mat{M} \partial_t \vec{u},
\end{equation}
where the second term on the right hand side has been estimated for the
semidiscretisation. Thus, the
fully discrete method is conservative across elements, if this is fulfilled by
the semidiscrete scheme. Considering stability,
\begin{equation}
\label{eq:stability-euler}
  \vec{u}_+^T \mat{M} \vec{u}_+
  = \vec{u}^T \mat{M} \vec{u}
    + 2 \Delta t \, \vec{u}^T \mat{M} \partial_t \vec{u}
    + (\Delta t)^2 \left(\partial_t \vec{u} \right)^T \mat{M} \partial_t \vec{u}.
\end{equation}
The second term is estimated in the same way as for the semidiscretisation.
However, the last term is non-negative, since $\mat{M}$ is positive definite.
Therefore, the fully discrete scheme may not yield the same estimates as the
continuous equation.

If possibility \ref{itm:operator-splitting}, i.e. the application of an operator
splitting method, is considered, a full time step using the explicit Euler
method is given by
\begin{equation}
\label{eq:filtered-euler}
  \vec{u}
  \mapsto \tilde{\vec{u}}_+ := \vec{u} + \Delta t \, \partial_t \vec{u}
  \mapsto \vec{u}_+ := \mat{F} \tilde{\vec{u}}_+.
\end{equation}
Therefore, equation \eqref{eq:stability-euler} holds for $\tilde{\vec{u}}_+$ instead
of $\vec{u}_+$. If the filter $\mat{F}$ reduces the norm of $\tilde{\vec{u}}_+$ by the
amount of the additional term $(\Delta t)^2 \left(\partial_t \vec{u} \right)^T
\mat{M} \partial_t \vec{u}$, the fully discrete scheme allows the same estimate
as the semidiscrete one.
Precisely, this idea is formulated in
\begin{lem}
\label{lem:operator-splitting}
  Rendering a conservative and stable semidiscretisation of the scalar conservation
  law \eqref{eq:scalar-CL}
  \begin{equation}
    \partial_t u + \partial_x f(u) = 0
  \end{equation}
  fully discrete by using an explicit Euler step with modal filtering
  \eqref{eq:filtered-euler} yields a conservative and stable scheme, if
  \begin{equation}
  \label{eq:filtered-euler-condition}
    \norm{\mat{F} \tilde{\vec{u}}_+}_M^2
    = \norm{\vec{u}}_M^2 + 2 \Delta t \scp{\vec{u}}{\partial_t \vec{u}}_M
    \leq \norm{ \tilde{\vec{u}}_+ }_M^2
    = \norm{\vec{u}}_M^2 + 2 \Delta t \scp{\vec{u}}{\partial_t \vec{u}}_M
      + (\Delta t)^2 \norm{\partial_t \vec{u}}_M^2.
  \end{equation}
  This condition can be fulfilled (per element) if
  \begin{itemize}
    \item 
    the rate of change $\partial_t \vec{u}$ is zero or
    
    \item
    the intermediate value $\tilde{\vec{u}}_+$ is not constant and the time step
    $\Delta t$ is small enough.
  \end{itemize}
\end{lem}
Of course, the same holds true if a strong-stability preserving (explicit) time
discretisation is chosen, since such a method can be written as a convex
combination of explicit Euler steps \citep{gottlieb2011strong}, and the filter
$\mat{F}$ is applied after each Euler step. However, these schemes will be
considered together with more general time discretisations in a forthcoming
article.

In order to fulfil condition \eqref{eq:filtered-euler-condition} of Lemma
\ref{lem:operator-splitting}, the filter strength $\epsilon$ (with time step 
$\Delta t$ included) can be adapted.
In a modal Legendre basis, the (exact) modal filter $\mat{F}$ may be written as
\begin{equation}
  \mat{F} = \diag{ \exp \left[ - \epsilon \, \lambda_n^s \right]_{n=0}^p },
\end{equation}
where $\lambda_n = \left(n (n+1) \right) \geq 0$. Representing the
polynomial given by $\tilde{\vec{u}}$ in a modal Legendre basis, i.e. as a
linear combination of Legendre polynomials $\phi_n$, condition
\eqref{eq:filtered-euler-condition} is
\begin{equation}
  \sum_{n=0}^p \exp[- 2 \epsilon \, \lambda_n^s] \, \tilde{u}_{+,n}^2 \, \norm{\phi_n}^2
  = RHS,
\end{equation}
where the right-hand side $\norm{u}_M^2 + 2 \Delta t \scp{\vec{u}}{\partial_t\vec{u}}$
is abbreviated as $RHS$. Using
\begin{equation}
  \exp[x] \geq 1 + x, \quad x \in \R,
\end{equation}
as a first order approximation, $\epsilon$ can be estimated by
\begin{equation}
\begin{aligned}
  & \sum_{n=0}^p (1 - 2 \epsilon \, \lambda_n^s) \, \tilde{u}_{+,n}^2 \, \norm{\phi_n}^2
    \leq RHS
  \\ \Leftrightarrow
  & \left( \sum_{n=0}^p \tilde{u}_{+,n}^2 \, \norm{\phi_n}^2 - RHS \right)
    \left( \sum_{n=0}^p 2 \lambda_n^s \, \tilde{u}_{+,n}^2 \, \norm{\phi_n}^2 \right)^{-1}
    \leq \epsilon,
\end{aligned}
\end{equation}
for $\sum_{n=0}^p 2 \lambda_n^s \, \tilde{u}_{+,n}^2 \, \norm{\phi_n}^2 > 0$.
Inserting 
\begin{equation}
\begin{aligned}
  \sum_{n=0}^p \tilde{u}_{+,n}^2 \, \norm{\phi_n}^2
  &=
  \norm{\vec{u}}_M^2 + 2 \Delta t \scp{\vec{u}}{\partial_t \vec{u}}_M
  + (\Delta t)^2 \norm{\partial_t \vec{u}}_M^2
  \\&=
  RHS + (\Delta t)^2 \norm{\partial_t \vec{u}}_M^2,
\end{aligned}
\end{equation}
this yields
\begin{lem}
\label{lem:estimate-epsilon}
  A necessary condition for the filter strength according to Lemma
  \ref{lem:operator-splitting} is
  \begin{equation}
  \label{eq:estimate-epsilon}
    \epsilon
    \geq
    (\Delta t)^2 \norm{\partial_t \vec{u}}_M^2
    \left( \sum_{n=0}^p 2 \lambda_n^s \, \tilde{u}_{+,n}^2 \, \norm{\phi_n}^2 \right)^{-1}.
  \end{equation}
\end{lem}

\subsection{Numerical experiments}

In this section, some numerical experiments will be performed to augment the
theoretical considerations.

\subsubsection{Linear advection with smooth solution}

The linear advection equation \eqref{eq:lin-adv}
\begin{equation}
  \partial_t u + \partial_x u = 0
  ,\quad
  u(0,x) = u_0(x) = \exp \left( -20 (x-1)^2 \right)
\end{equation}
is solved by an SBP CPR semidiscretisation \eqref{eq:lin-adv-semidisc} using
$N = 8$ elements with nodal Gauß-Legendre bases of degree $\leq p = 7$ in the
domain $[0, 2]$ with periodic boundary conditions. An explicit Euler method
using $12 \cdot 10^4$ time steps is used to advance the solution in the time
interval $[0, 10]$ and a central numerical flux $\fnum(u_-,u_+) = (u_- + u_+) / 2$
is applied.

\begin{figure}[!ht]
  \centering
  \begin{subfigure}[b]{0.49\textwidth}
    \includegraphics[width=\textwidth]{%
      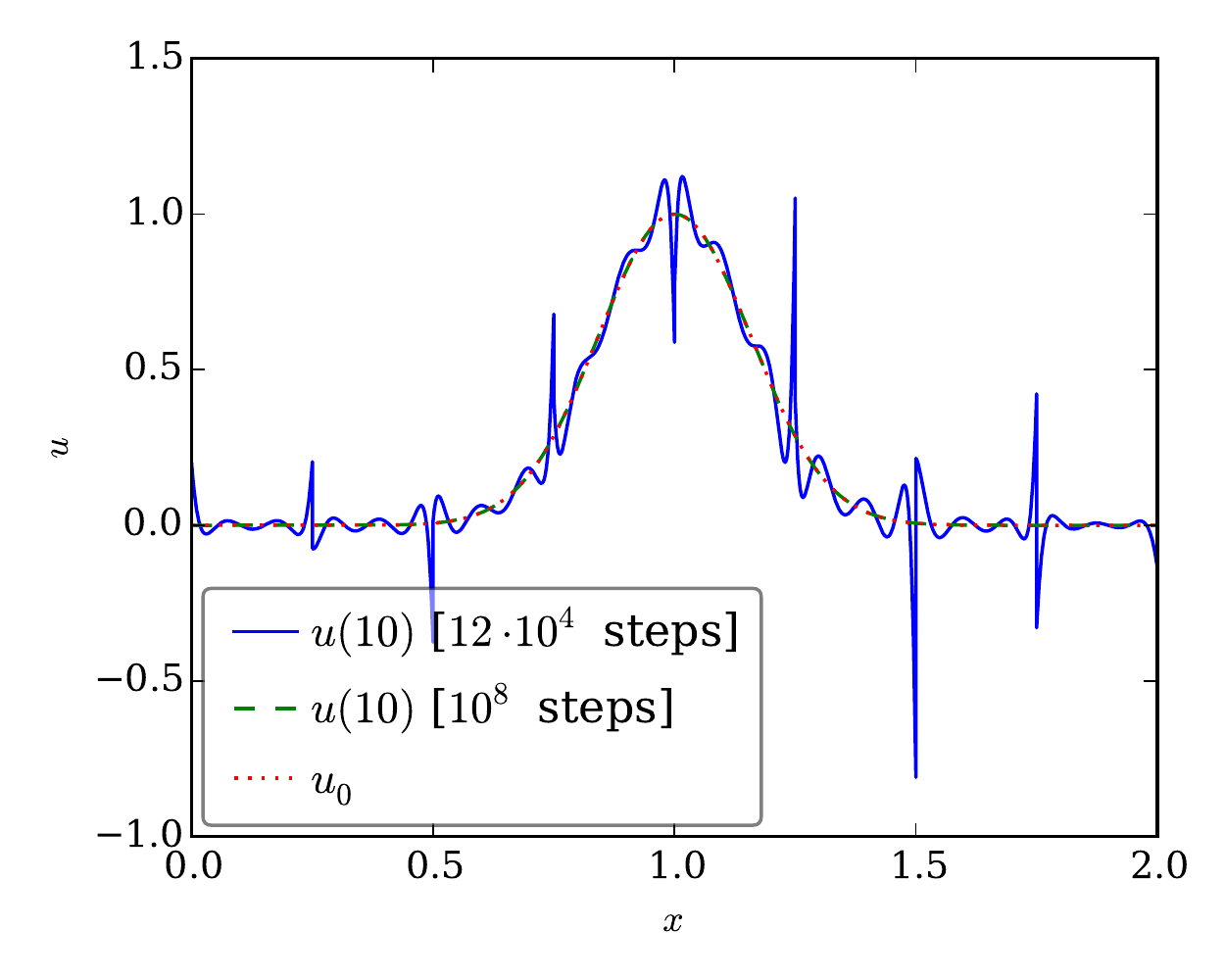}
    \caption{Solution computed without modal filtering.}
    \label{fig:linear_gauss_no_filtering_solution}
  \end{subfigure}%
  ~
  \begin{subfigure}[b]{0.49\textwidth}
    \includegraphics[width=\textwidth]{%
      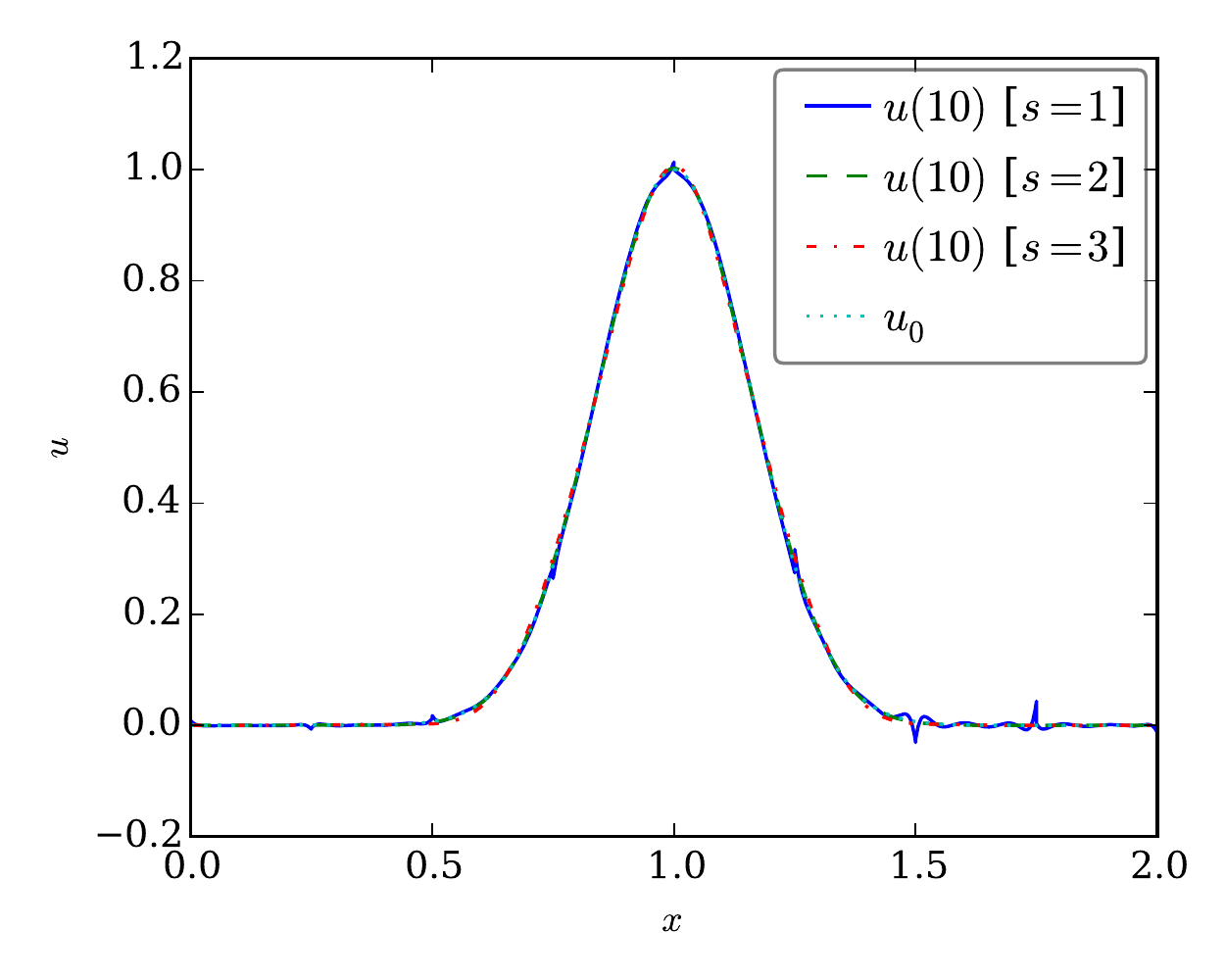}
    \caption{Solution computed using adaptive modal filtering.}
    \label{fig:linear_gauss_mf_adaptive_solution}
  \end{subfigure}%
  \\
  \begin{subfigure}[b]{0.49\textwidth}
    \includegraphics[width=\textwidth]{%
      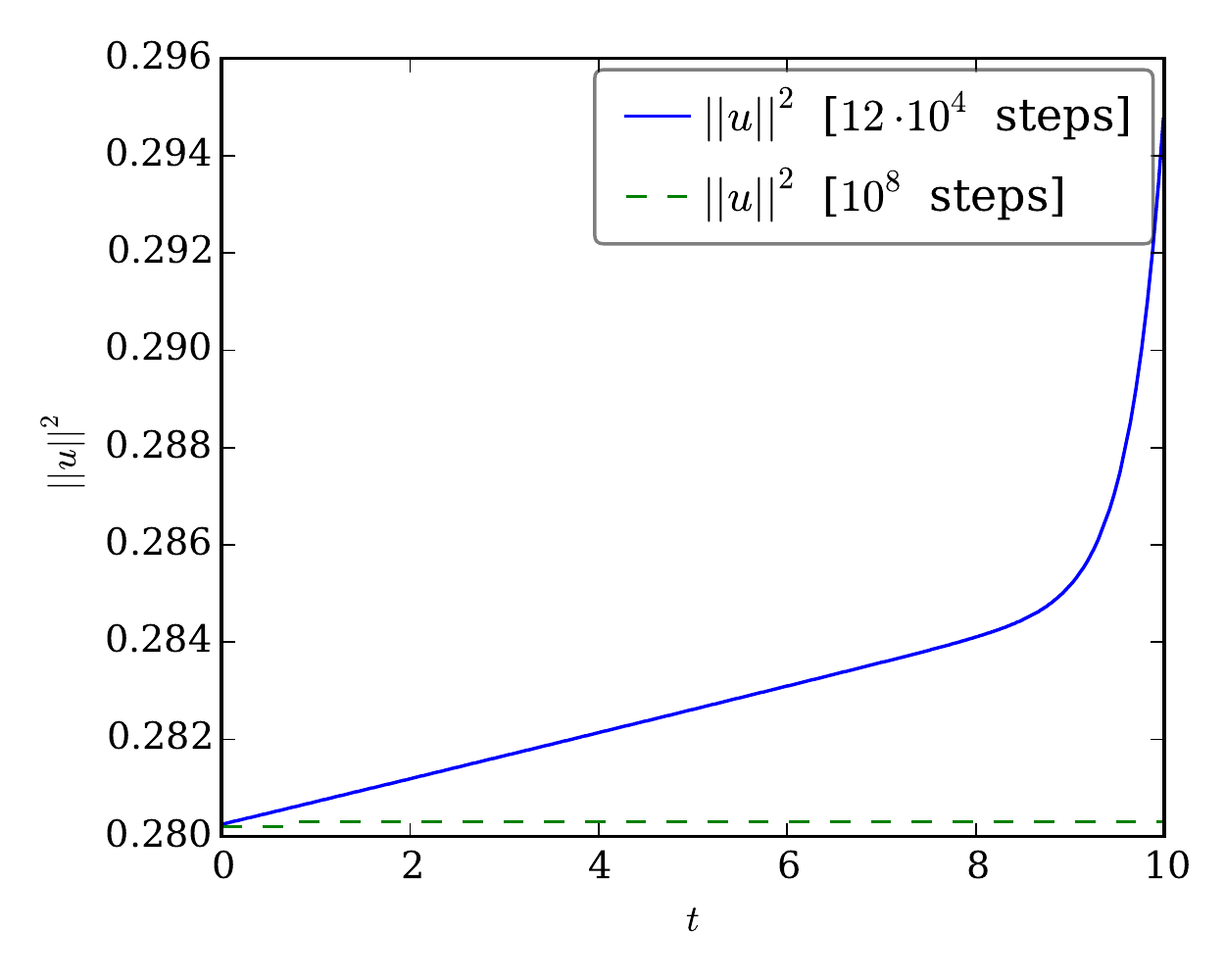}
    \caption{Energy of the solution computed without additional modal filtering.}
    \label{fig:linear_gauss_no_filtering_energy}
  \end{subfigure}%
  ~
  \begin{subfigure}[b]{0.49\textwidth}
    \includegraphics[width=\textwidth]{%
      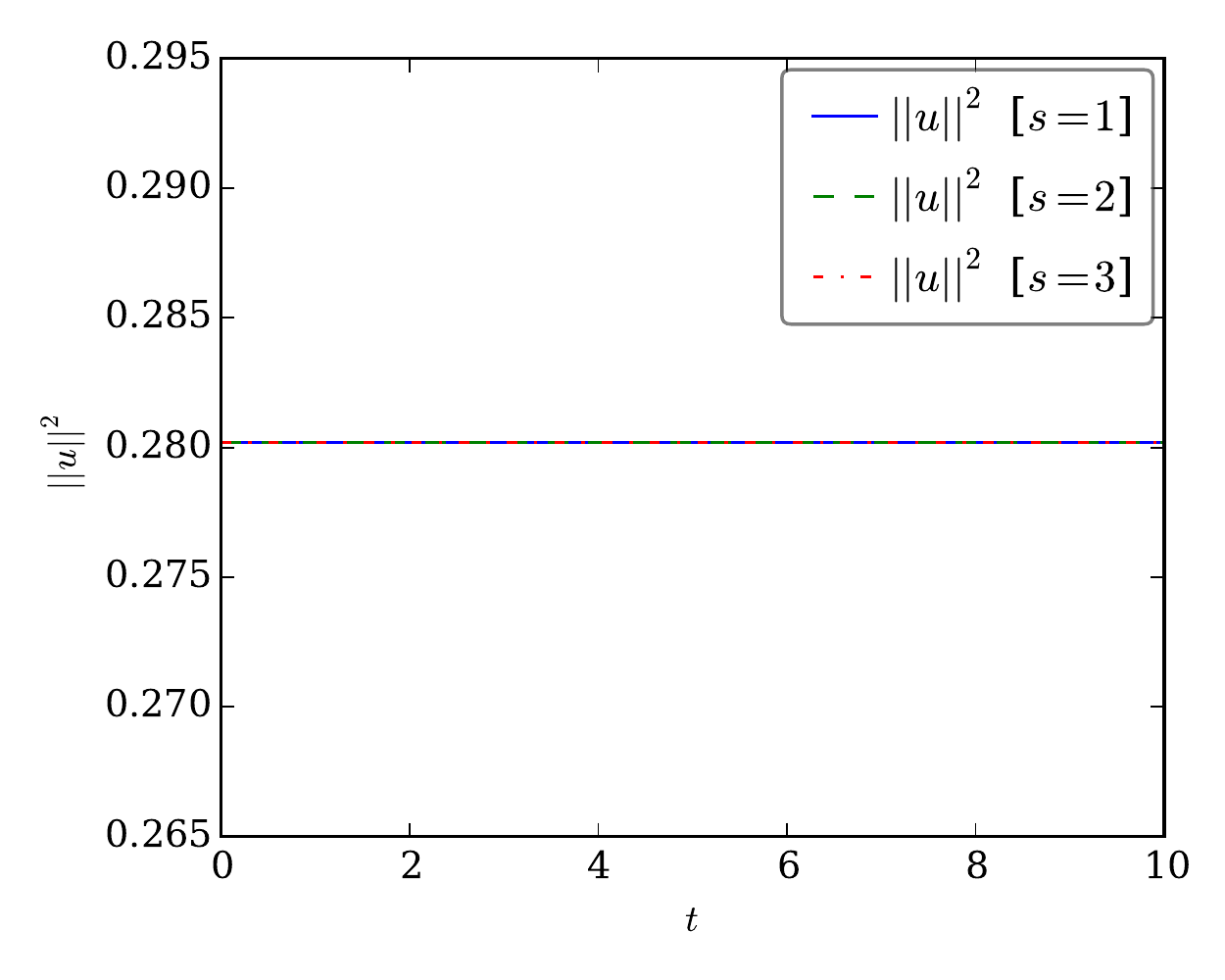}
    \caption{Energy of the solution computed using adaptive modal filtering.}
    \label{fig:linear_gauss_mf_adaptive_energy}
  \end{subfigure}%
  \caption{Numerical results for linear advection using $N = 8$ elements with
           polynomials of degree $\leq p = 7$.
           On the left hand side, no modal filtering has been used, whereas
           adaptive  modal filters of orders $s=1$, $s=2$, and $s=3$ have been
           applied for the right hand side.
           }
  \label{fig:linear_gauss}
\end{figure}

Figure \ref{fig:linear_gauss_no_filtering_solution} shows the oscillating solution
without modal filtering using $12 \cdot 10^4$ time steps. The corresponding
increasing energy is plotted in Figure \ref{fig:linear_gauss_no_filtering_energy}.
However, increasing the number of time steps to $10^6$ drastically reduces both
the increase of energy and the oscillations.

On the other hand, applying modal filters of order $s \in \set{1,2,3}$ with
adaptively chosen strength $\epsilon$ according to Lemma \ref{lem:estimate-epsilon}
results in constant energy in Figure \ref{fig:linear_gauss_mf_adaptive_energy}
and non-oscillatory solutions in Figure \ref{fig:linear_gauss_mf_adaptive_solution}.
Only the numerical solution computed with filters of order $s = 1$ shows
slight peaks in the smooth part.

\begin{figure}[!hp]
  \centering
  \begin{subfigure}[b]{0.49\textwidth}
    \includegraphics[width=\textwidth]{%
      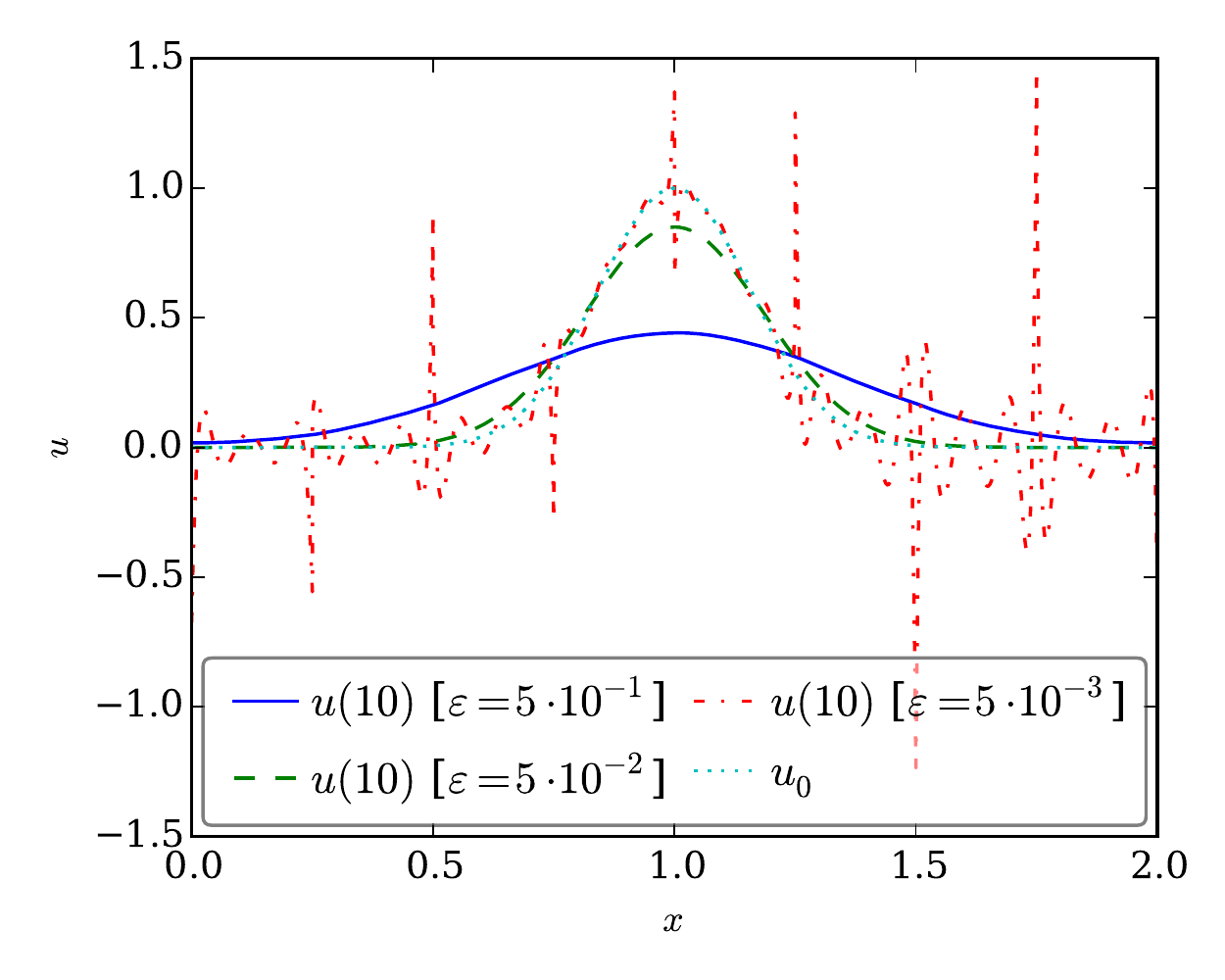}
    \caption{Solutions for order $s = 1$.}
  \end{subfigure}%
  ~
  \begin{subfigure}[b]{0.49\textwidth}
    \includegraphics[width=\textwidth]{%
      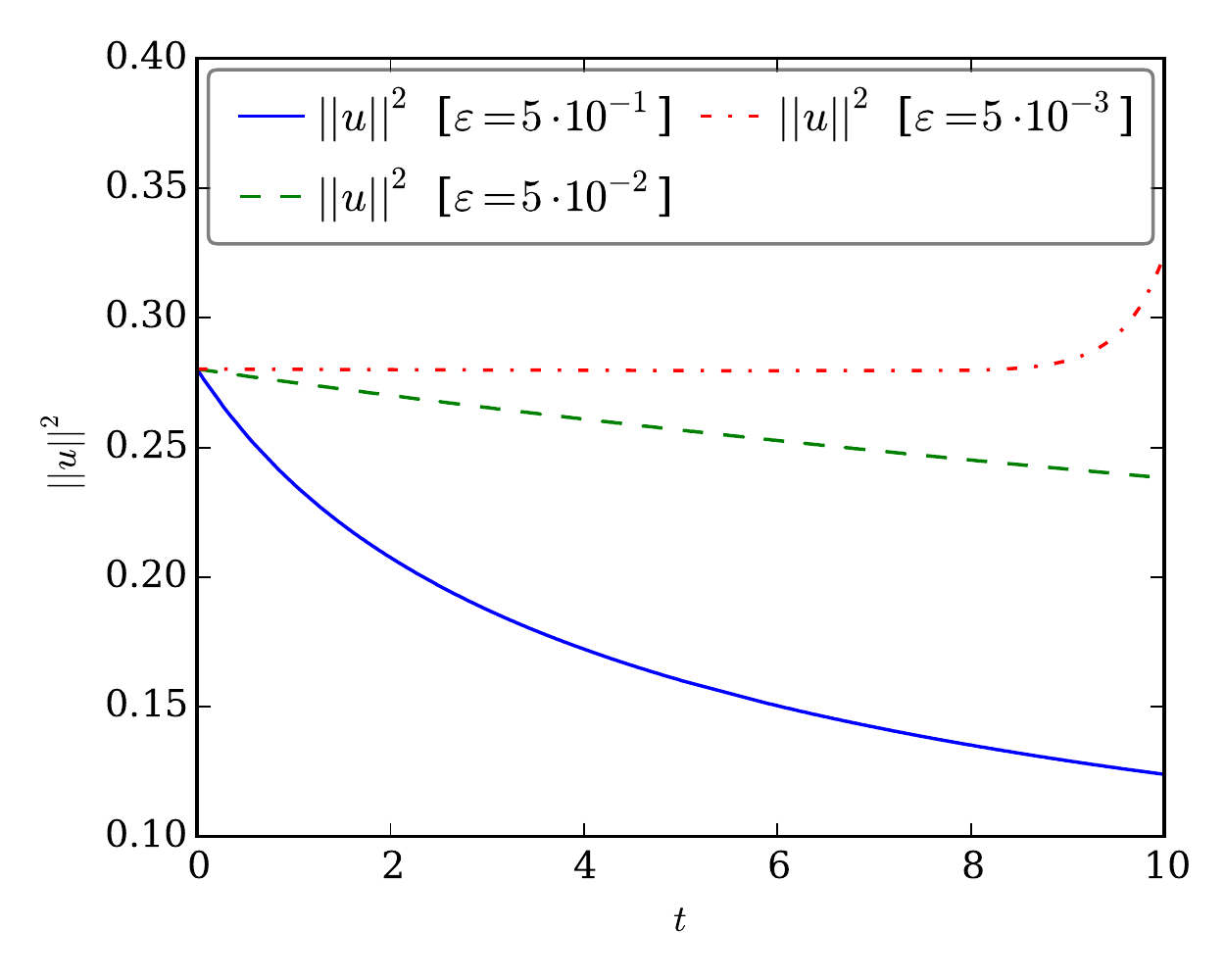}
    \caption{Energies for order $s = 1$.}
  \end{subfigure}%
  \\
  \begin{subfigure}[b]{0.49\textwidth}
    \includegraphics[width=\textwidth]{%
      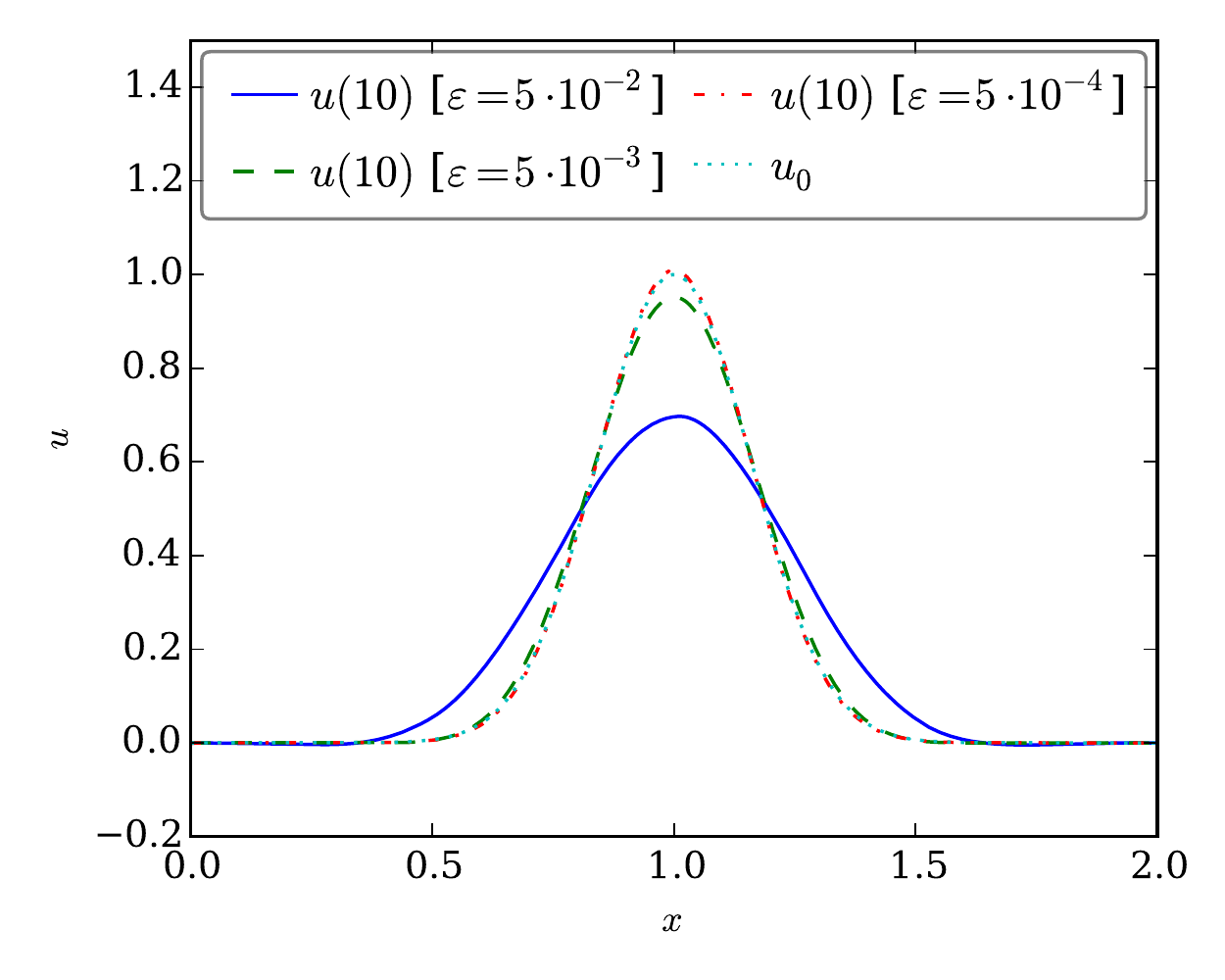}
    \caption{Solutions for order $s = 2$.}
  \end{subfigure}%
  ~
  \begin{subfigure}[b]{0.49\textwidth}
    \includegraphics[width=\textwidth]{%
      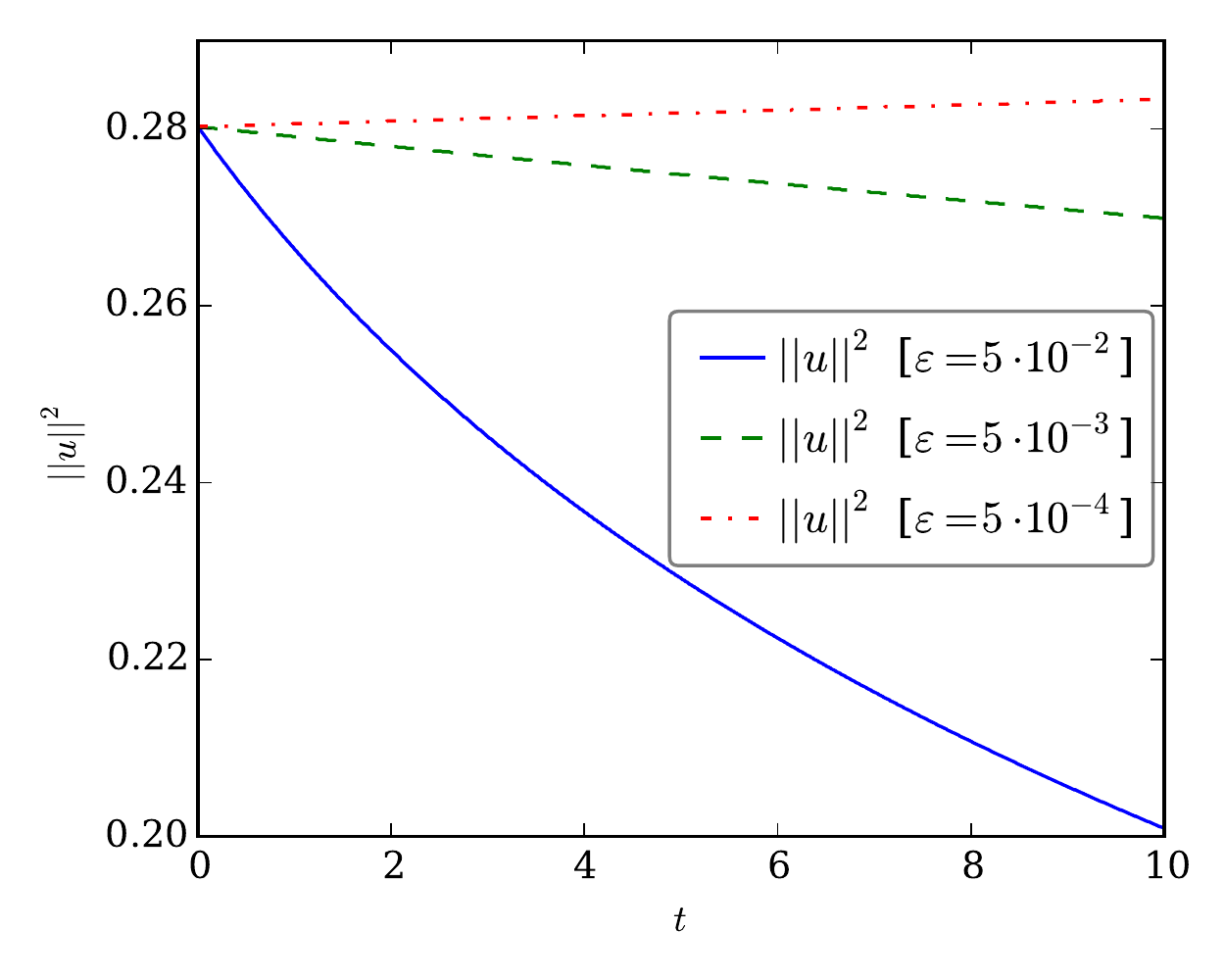}
    \caption{Energies for order $s = 2$.}
  \end{subfigure}%
  \\
  \begin{subfigure}[b]{0.49\textwidth}
    \includegraphics[width=\textwidth]{%
      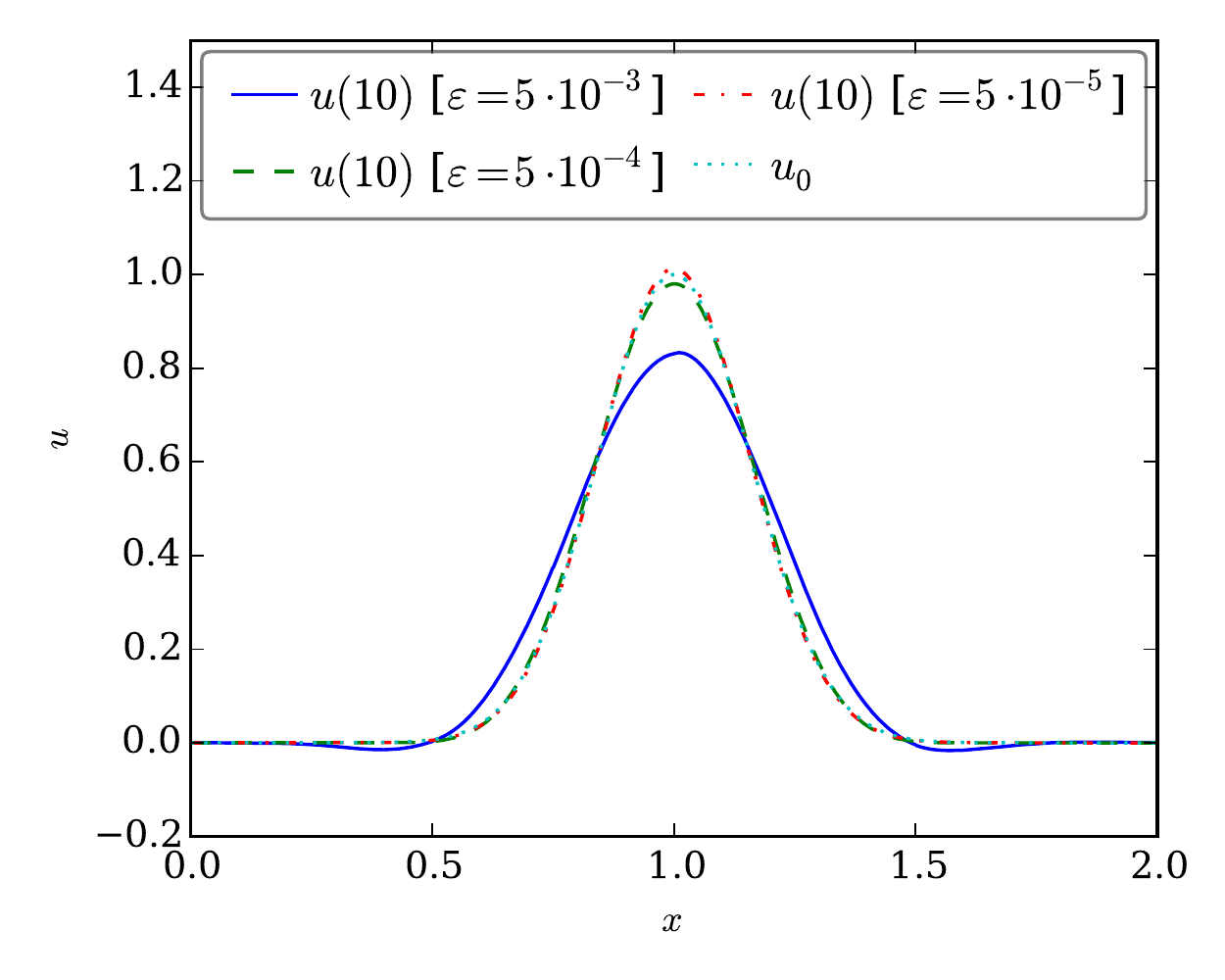}
    \caption{Solutions for order $s = 3$.}
  \end{subfigure}%
  ~
  \begin{subfigure}[b]{0.49\textwidth}
    \includegraphics[width=\textwidth]{%
      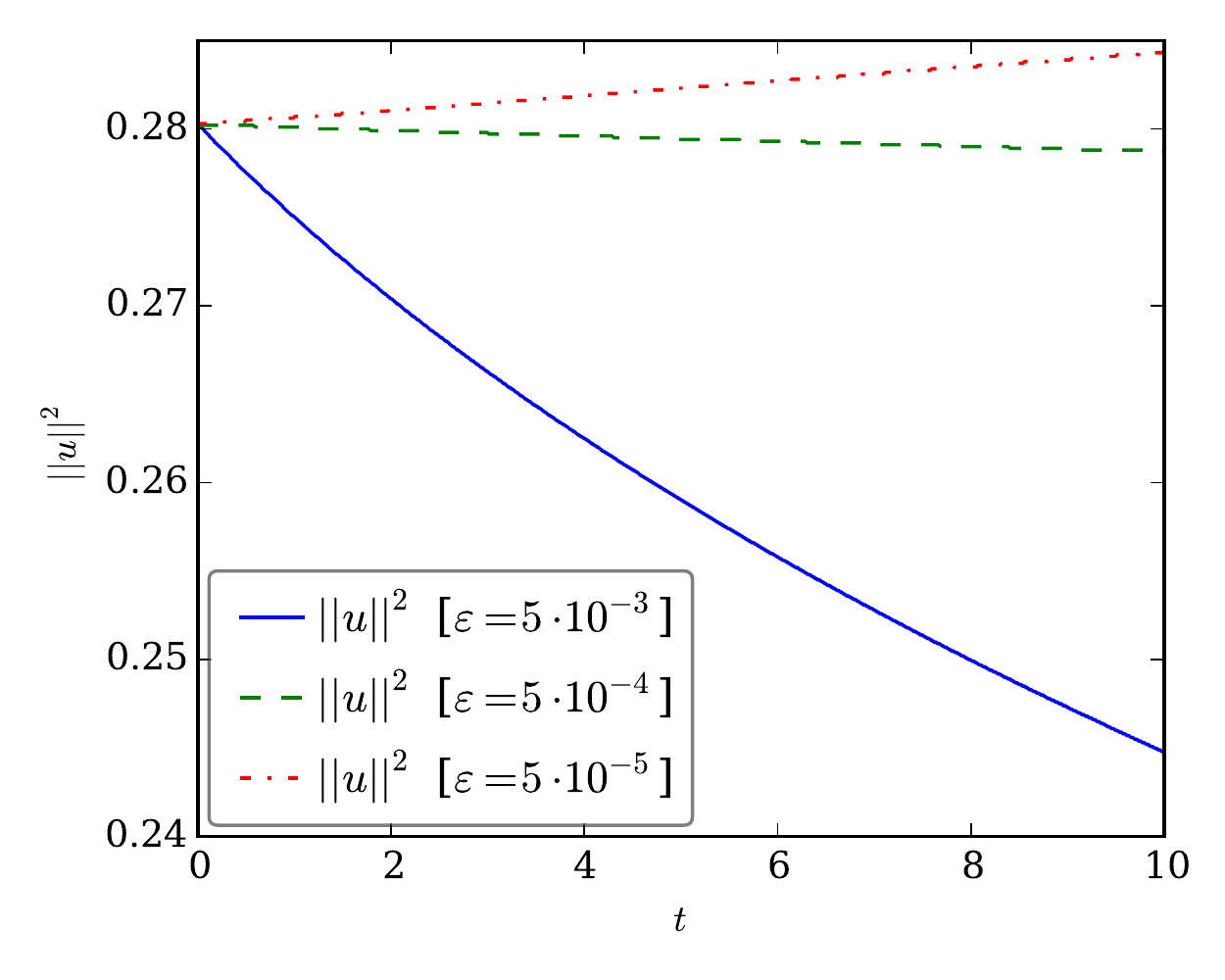}
    \caption{Energies for order $s = 3$.}
  \end{subfigure}%
  \caption{Numerical results for linear advection using $N = 8$ elements with
           polynomials of degree $\leq p = 7$ and modal filtering of orders
           $s \in \set{1,2,3}$ with various strengths $\epsilon$.
           On the left hand side, the solutions $u$ are shown, accompanied by the
           corresponding energies $\norm{u}^2$ on the right hand side.
           }
  \label{fig:linear_gauss_comparison}
\end{figure}

A simple application of modal filtering with constant order $s$ and strength
$\epsilon$ has a stabilising effect, as can be seen in Figure
\ref{fig:linear_gauss_comparison}. The dissipation increases with increasing
strength $\epsilon$ and order $s$, respectively. However, numerous experiments
and fine tuning of the parameters by hand is required in order to get
acceptable results. Therefore, the adaptive choice proposed in this work is
advantageous.

\subsubsection{Linear advection with discontinuous solution}

The influence of the presence of discontinuities is investigated using the linear
advection equation \eqref{eq:lin-adv}
\begin{equation}
\label{eq:lin-adv-jump}
  \partial_t u + \partial_x u = 0
  ,\qquad
  u(0, x) = u_0(x) =
  \begin{cases}
    1, & x \in [0.5, 1],\\
    0, & \text{otherwise},
  \end{cases}
\end{equation}
and a semidiscretisation \eqref{eq:lin-adv-semidisc} using nodal Gauß-Legendre
bases of degree $\leq p = 7$ on $N = 8$ elements with an upwind numerical flux
$\fnum(u_-,u_+) = u_-$. The domain $[0,2]$ is equipped with periodic boundary
conditions and the solution is advanced in time $t \in [0, 8]$ by
an explicit Euler method.

\begin{figure}[!htb]
  \centering
  \begin{subfigure}[b]{0.49\textwidth}
    \includegraphics[width=\textwidth]{%
      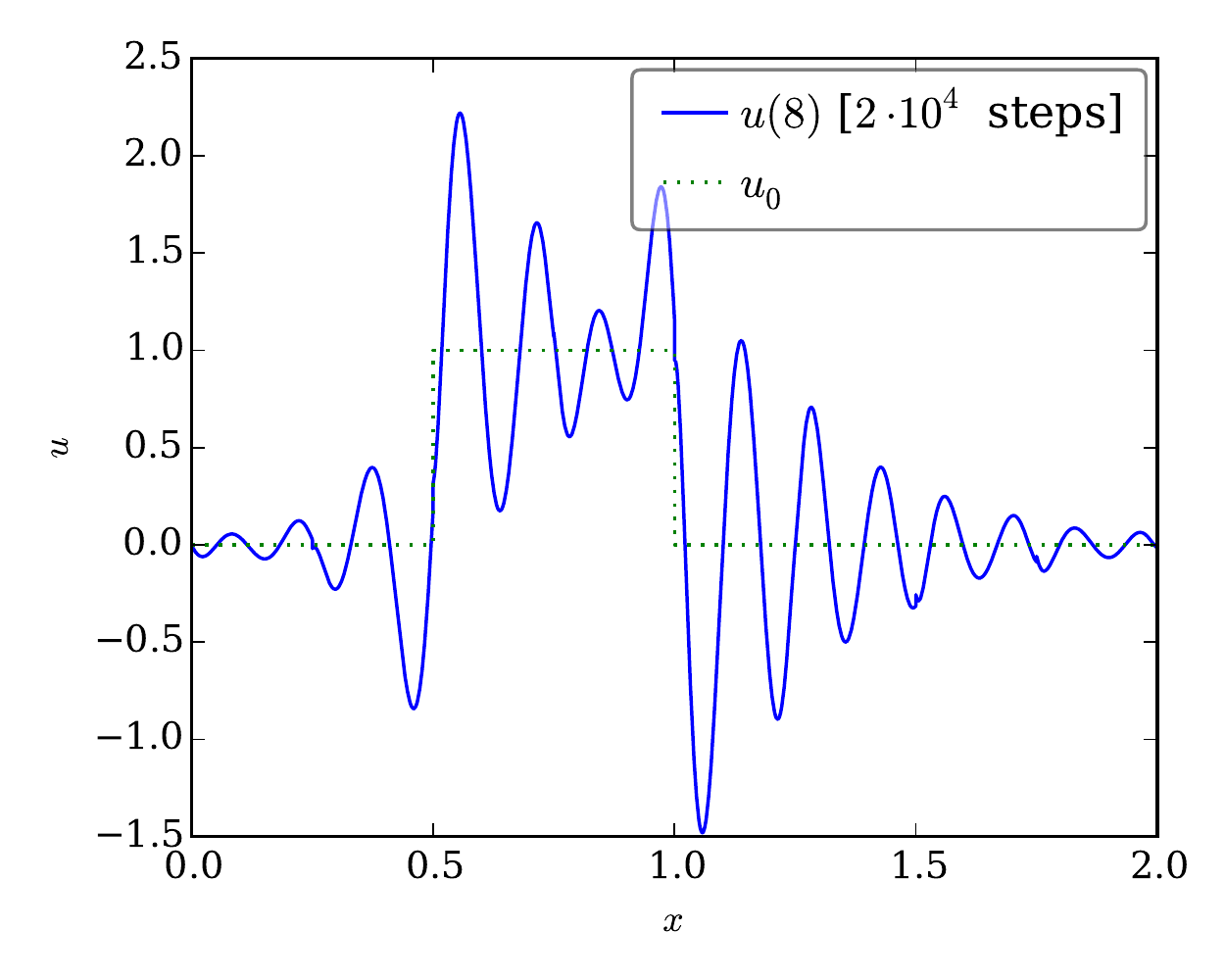}
    \caption{Solution computed without modal filtering.}
    \label{fig:linear_jump_no_filtering_solution}
  \end{subfigure}%
  ~
  \begin{subfigure}[b]{0.49\textwidth}
    \includegraphics[width=\textwidth]{%
      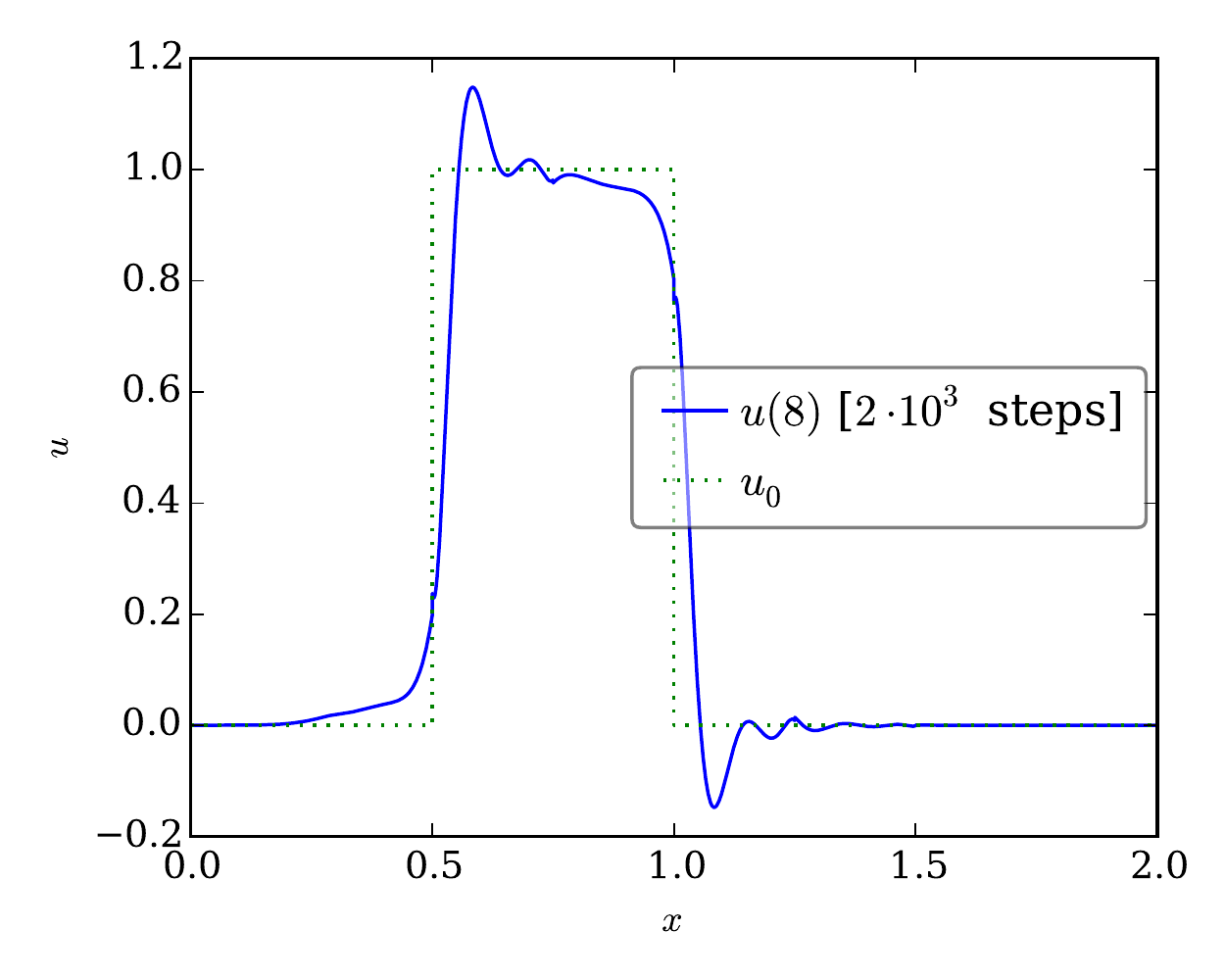}
    \caption{Solution computed using adaptive modal filtering.}
    \label{fig:linear_jump_ad_adaptive_solution}
  \end{subfigure}%
  \\
  \begin{subfigure}[b]{0.49\textwidth}
    \includegraphics[width=\textwidth]{%
      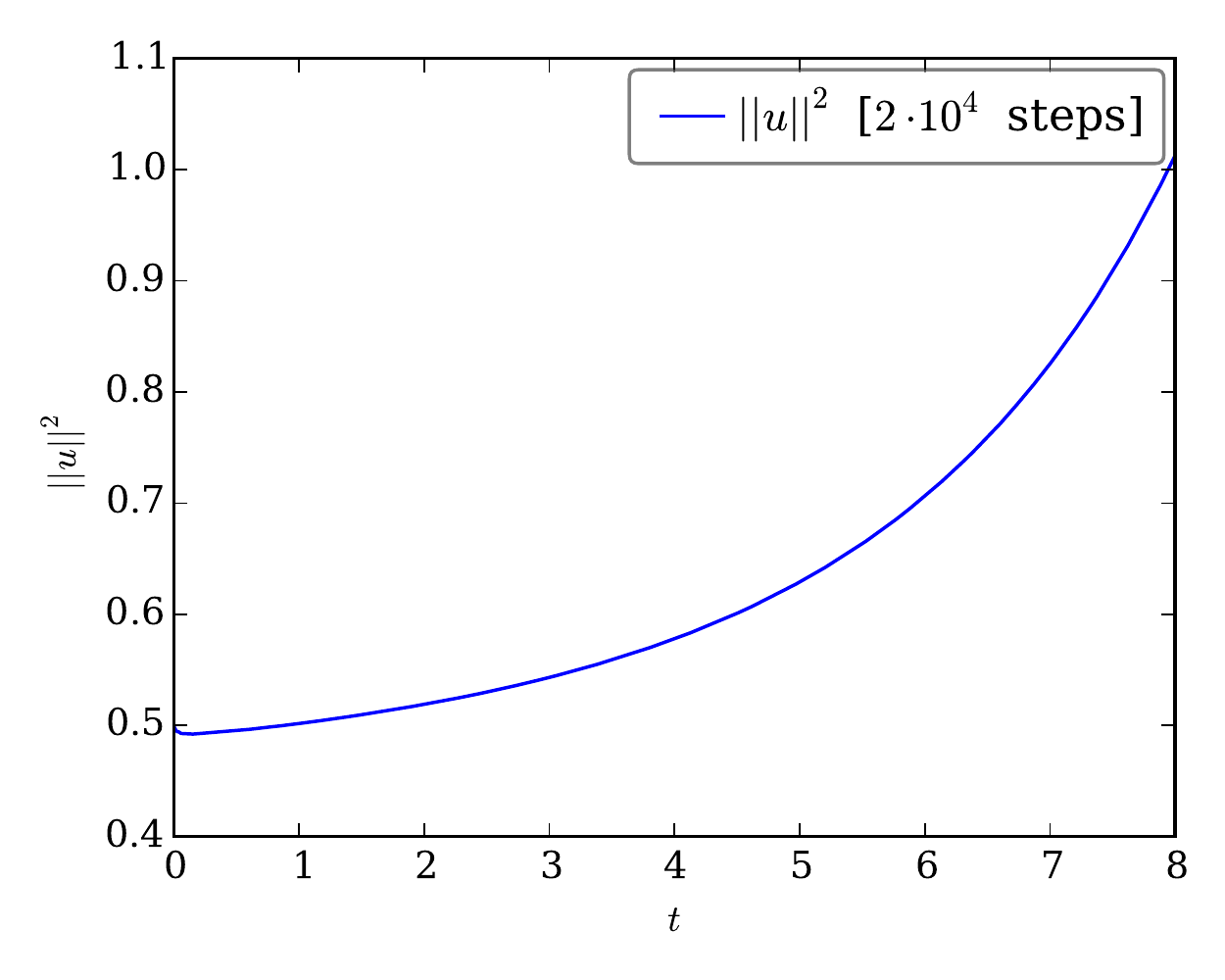}
    \caption{Energy of the solution computed without modal filtering.}
    \label{fig:linear_jump_no_filtering_energy}
  \end{subfigure}%
  ~
  \begin{subfigure}[b]{0.49\textwidth}
    \includegraphics[width=\textwidth]{%
      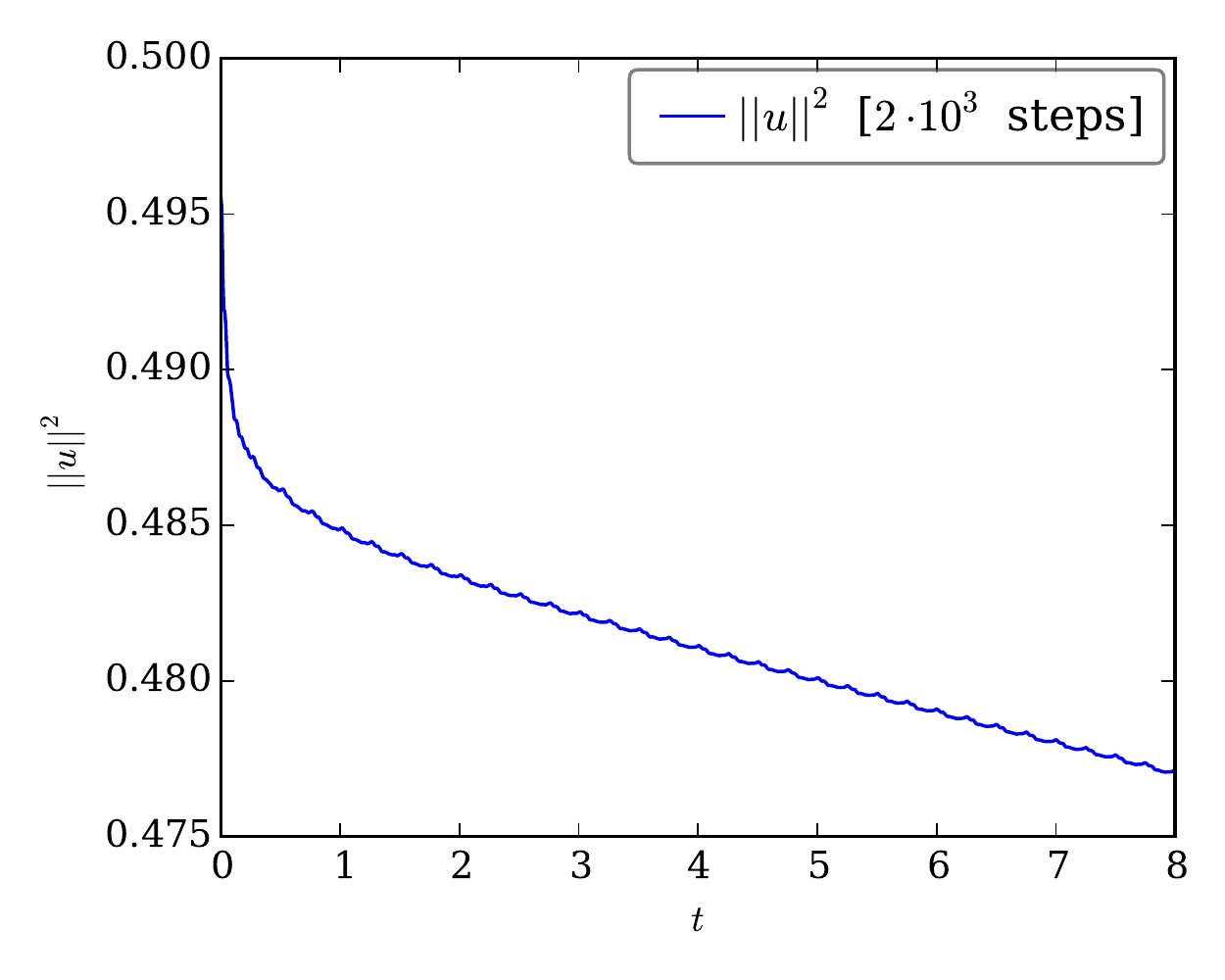}
    \caption{Energy of the solution computed using adaptive modal filtering.}
    \label{fig:linear_jump_ad_adaptive_energy}
  \end{subfigure}%
  \caption{Numerical results for linear advection using $N = 16$ elements with
           polynomials of degree $\leq p = 15$.
           On the left hand side, no modal filtering has been used, whereas adaptive
           modal filters of order $s=1$ have been used for the right hand side.
           }
  \label{fig:linear_jump}
\end{figure}

Using $2 \cdot 10^4$ time steps, the solution without modal filtering in Figure
\ref{fig:linear_jump_no_filtering_solution} is oscillatory with increasing
energy (Figure \ref{fig:linear_jump_no_filtering_energy}). However, applying
modal filtering with adaptively chosen strength $\epsilon$ yields a much less
oscillatory solution in Figure \ref{fig:linear_jump_ad_adaptive_solution}
with slightly decreasing energy in Figure \ref{fig:linear_jump_ad_adaptive_energy}
using only $2 \cdot 10^3$ time steps.

\subsubsection{Burgers' equation}

A nonlinear model problem is given by Burgers' equation \eqref{eq:Burgers}
\begin{equation}
  \partial_t u + \partial_x \frac{u^2}{2} = 0
  ,\quad
  u(0,x) = u_0(x) = \sin \pi x + 0.01
\end{equation}
in the domain $[0, 2]$ with periodic boundary conditions. The semidiscretisation
\eqref{eq:Burgers-semidisc} with $N = 16$ elements using Gauß-Legendre bases of
degree $\leq p = 15$ and a local Lax-Friedrichs numerical flux
$\fnum(u_-,u_+) = \frac{u_-^2 + u_+^2}{4}
- \frac{ \max \set{ \abs{u_-}, \abs{u_+} } }{2} (u_+ - u_-)$ is applied together
with an explicit Euler method.

\begin{figure}[!htb]
  \centering
  \begin{subfigure}[b]{0.49\textwidth}
    \includegraphics[width=\textwidth]{%
      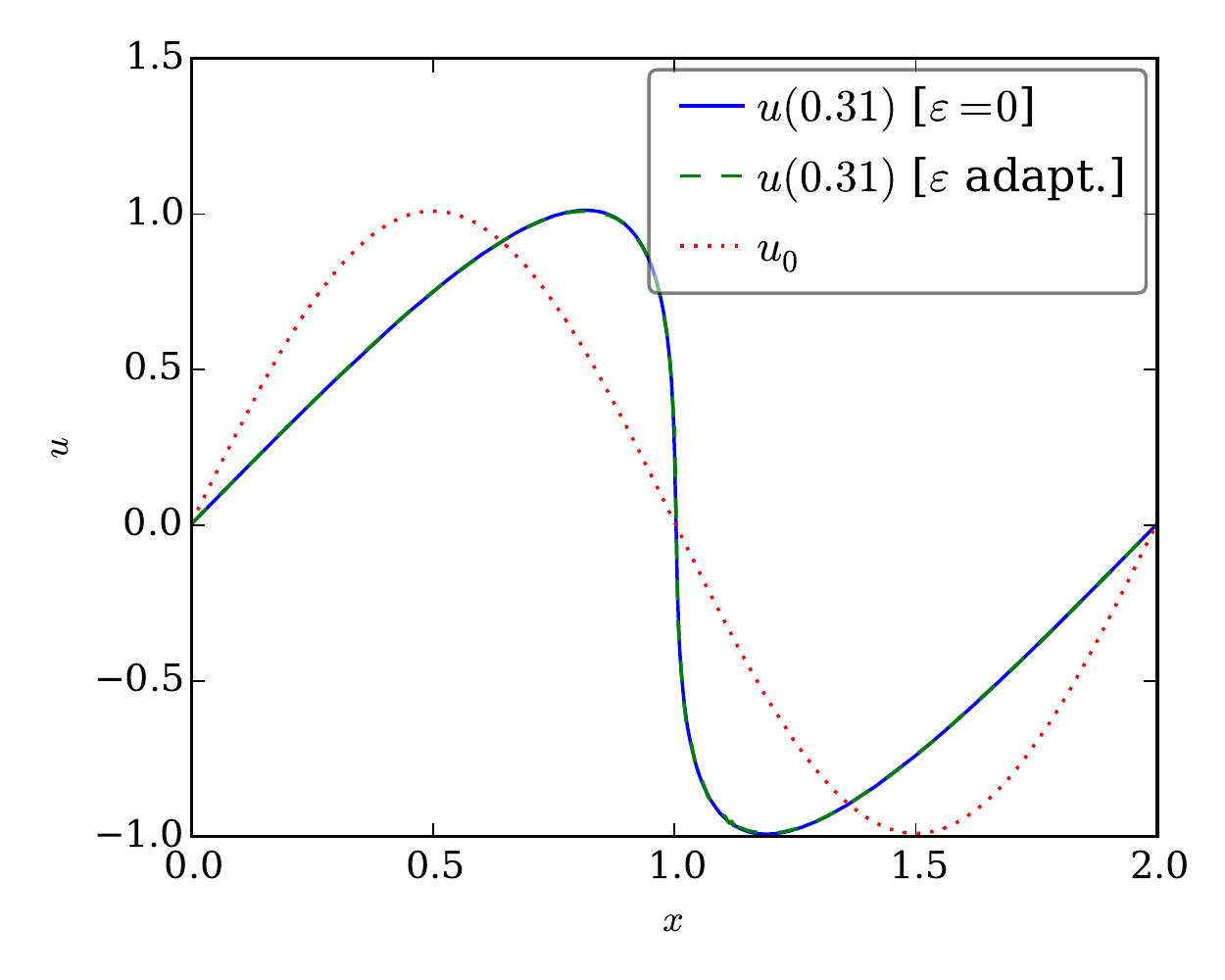}
    \caption{Solution at time $t = 0.31$.}
    \label{fig:burgers_sin_smooth_solution}
  \end{subfigure}%
  ~
  \begin{subfigure}[b]{0.49\textwidth}
    \includegraphics[width=\textwidth]{%
      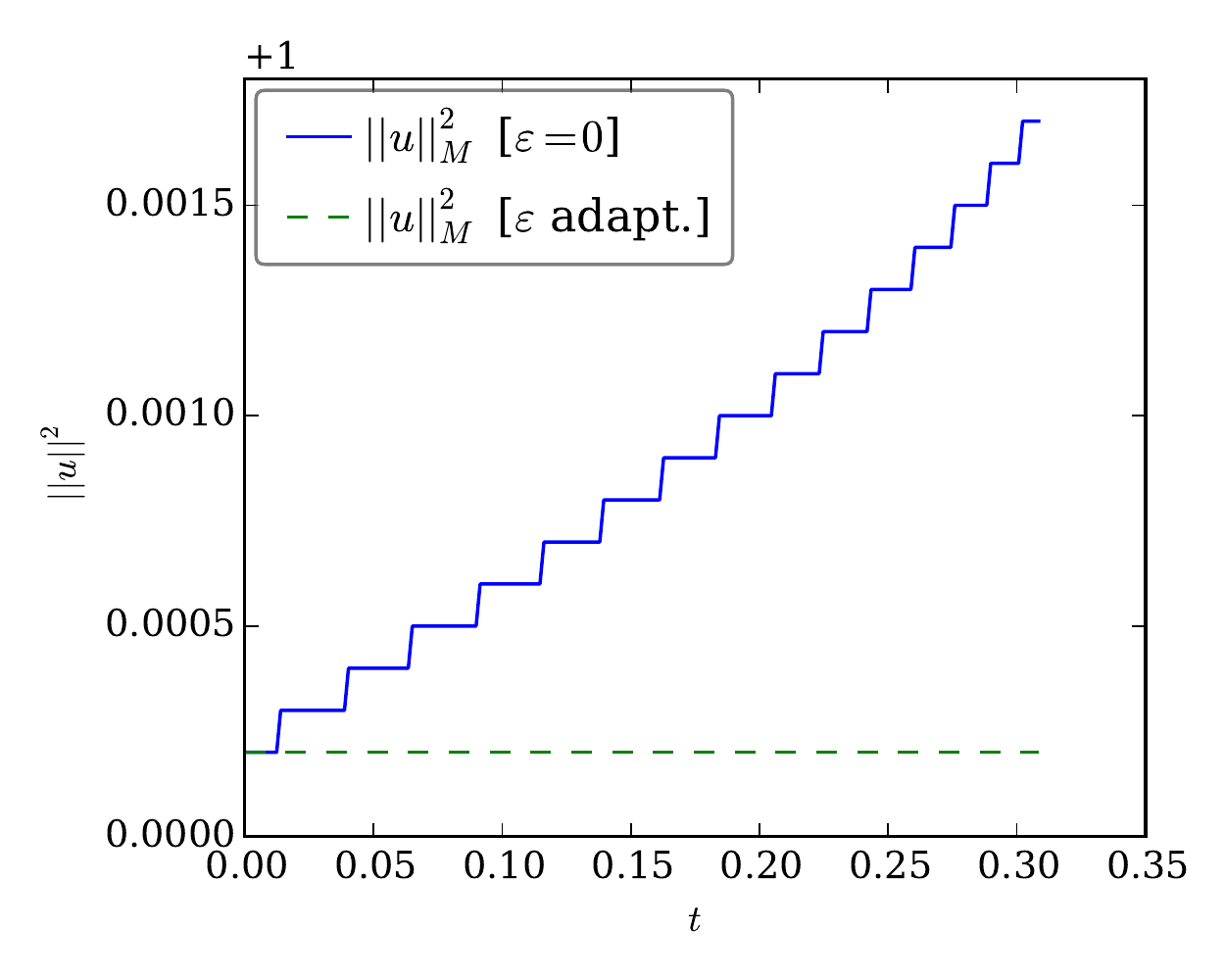}
    \caption{Energy in the time interval $[0, 0.31]$.}
    \label{fig:burgers_sin_smooth_energy}
  \end{subfigure}%
  \\
  \begin{subfigure}[b]{0.49\textwidth}
    \includegraphics[width=\textwidth]{%
      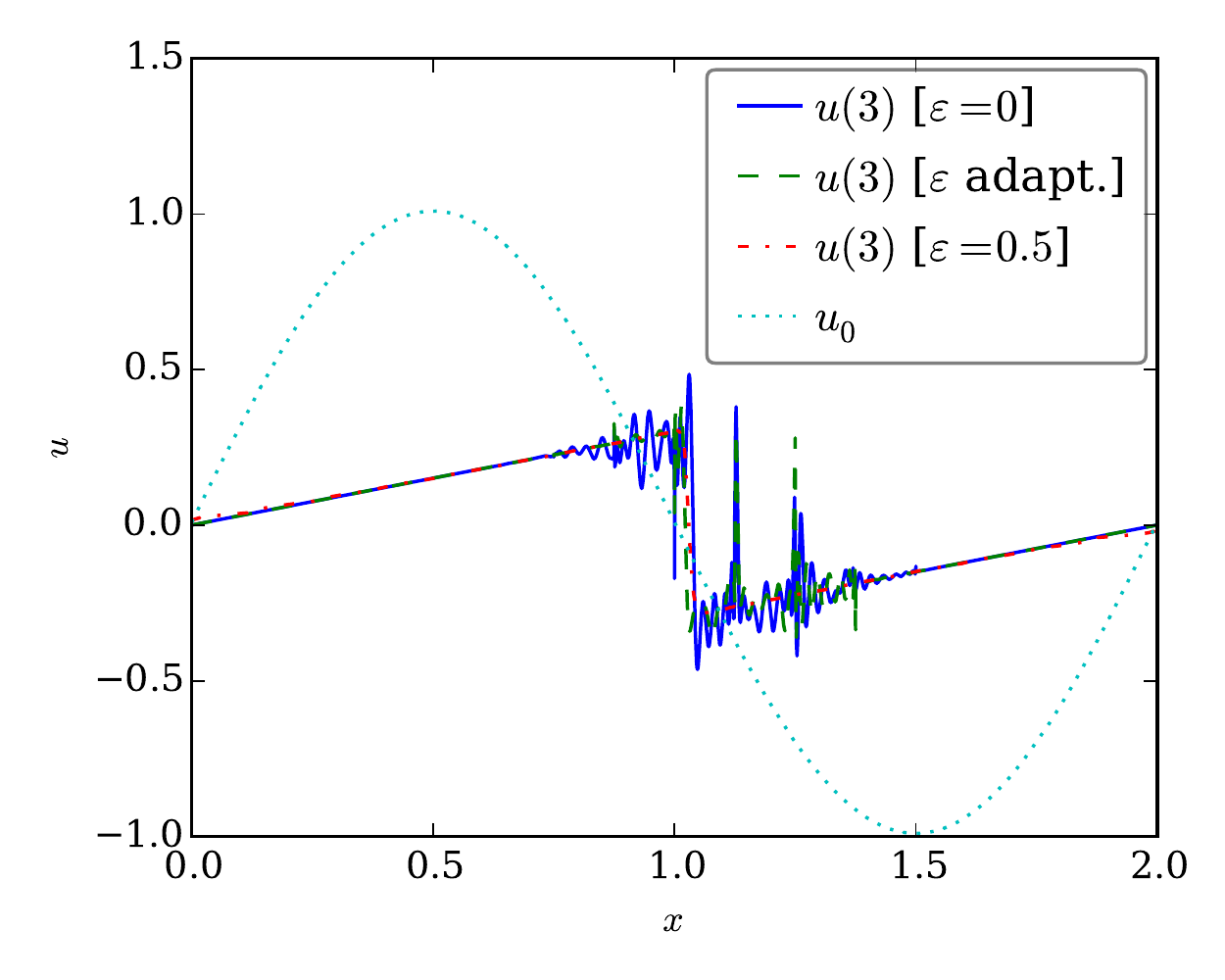}
    \caption{Solution at time $t = 3$.}
    \label{fig:burgers_sin_rough_solution}
  \end{subfigure}%
  ~
  \begin{subfigure}[b]{0.49\textwidth}
    \includegraphics[width=\textwidth]{%
      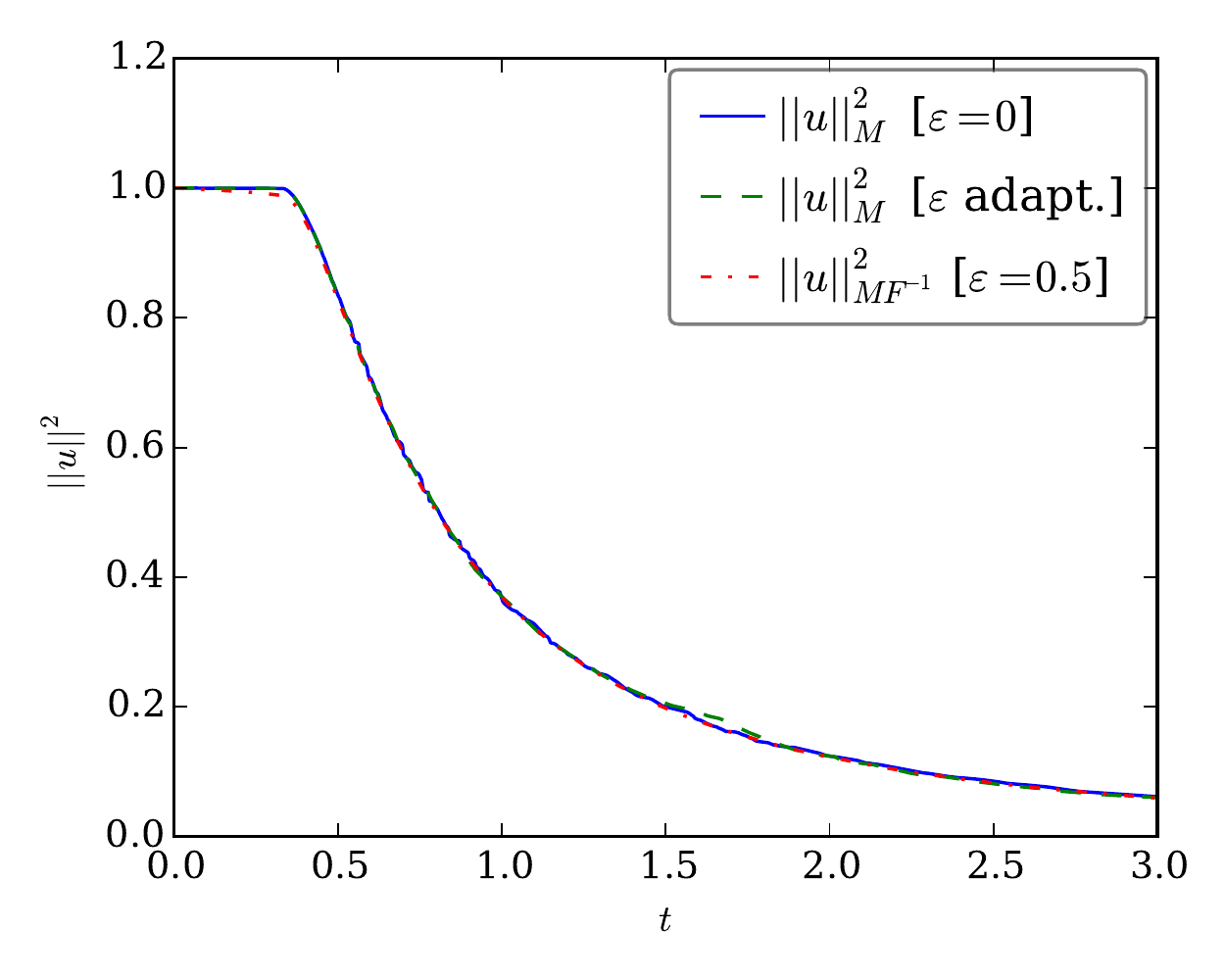}
    \caption{Energy in the time interval $[0, 3]$.}
    \label{fig:burgers_sin_rough_energy}
  \end{subfigure}%
  \caption{Numerical results for Burgers' equation using $N = 16$ elements with
           polynomials of degree $\leq p = 15$.
           The solutions and the energy are plotted on the left hand side and
           right hand side, respectively, for the time intervals $[0, 0.31]$ (first
           row) and $[0, 3]$ (second row).
           }
  \label{fig:burgers_sin}
\end{figure}

The numerical solution up to $t = 0.31$ is computed using $200$ time steps and
plotted in Figure \ref{fig:burgers_sin_smooth_solution} with and without
modal filtering. Both solutions coincide visually but the energy increases without
modal filtering. However, application of adaptive modal filtering yields a
constant energy, as expected.

Contrary, both the solutions and the energy on the complete time interval $[0, 3]$
using $15 \cdot 10^3$ steps are visually indistinguishable, as can be seen in
Figures \ref{fig:burgers_sin_rough_solution} and \ref{fig:burgers_sin_rough_energy}.
Thus, the adaptive modal filtering is not able to remove the oscillations triggered
by the developing discontinuity. However, a simple application of a modal filter
with fixed strength $\epsilon = 0.5$ results in a non-oscillatory solution.

  \section{Filtering the time derivative}
\label{sec:filter-derivative}

In this section, possibility \ref{itm:filter-derivative} of section
\ref{sec:discrete-setting} is discussed, i.e. the application of a modal
filter $\mat{F}$ to the time derivative of $\vec{u}$.

\subsection{Semidiscrete and discrete estimates}

In the following, a conservative and stable semidiscretisation of the scalar conservation law
\begin{equation}
  \partial_t u + \partial_x f(u) = 0
\end{equation}
relying on elements is assumed, e.g. a CPR method using SBP operators for the
linear advection equation with constant speed or Burgers' equation
\citep{ranocha2016summation}. This yields on each element an ordinary
differential equation
\begin{equation}
  \od{}{t} \vec{u} = g(\vec{u}).
\end{equation}
Applying the filter $\mat{F}$ to the time derivative results in
\begin{equation}
  \od{}{t} \vec{u} = \mat{F} g(\vec{u}).
\end{equation}

If constant functions are eigenvectors of the filter with eigenvalue $1$, i.e.
the filter does not change constants, and the filter is self-adjoint with respect
to the mass matrix $\mat{M}$, i.e. $\mat{M} \mat{F} = \mat{F}[^T] \mat{M}$, then the rate of
change of the integral of $\vec{u}$ over an element is
\begin{equation}
  \od{}{t} \vec{1}^T \mat{M} \vec{u}
  = \vec{1}^T \mat{M} \od{}{t} \vec{u}
  = \vec{1}^T \mat{M} \mat{F} g(\vec{u})
  = \vec{1}^T \mat{F}[^T] \mat{M} g(\vec{u})
  = \vec{1}^T \mat{M} g(\vec{u}),
\end{equation}
which is the same as for the semidiscretisation without filtering.

Investigating stability in the norm $\norm{\cdot}_M$ induced by $\mat{M}$ would
start canonically with
\begin{equation}
  \frac{1}{2} \od{}{t} \norm{\vec{u}}_M^2
  = \vec{u}^T \mat{M} \od{}{t} \vec{u}
  = \vec{u}^T \mat{M} \mat{F} g(\vec{u}).
\end{equation}
However, there do not seem to be simple estimates for this rate of change.
Contrary, the rate of change of the norm induced by $\mat{M} \mat{F}[^{-1}]$ (if
$\mat{F}$ is invertible and $\mat{M} \mat{F}[^{-1}]$ induces a scalar product) can be easily
estimated as
\begin{equation}
  \frac{1}{2} \od{}{t} \norm{\vec{u}}_{M F^{-1}}^2
  = \vec{u}^T \mat{M} \mat{F}[^{-1}] \od{}{t} \vec{u}
  = \vec{u}^T \mat{M} g(\vec{u}),
\end{equation}
i.e. the same estimate as for the semidiscretisation without filtering.

Assuming $\mat{F}$ is invertible, the bilinear form induced by $\mat{M} \mat{F}[^{-1}]$ is
symmetric iff for all $\vec{v}, \vec{u}$
\begin{equation}
  \vec{v}^T \mat{M} \mat{F} \vec{u} = \vec{u}^T \mat{M} \mat{F} \vec{v},
\end{equation}
i.e. $\mat{F}$ is $\mat{M}$-self-adjoint $\mat{M} \mat{F} = \mat{F}[^T] \mat{M}$.
Assuming $\mat{M}$ and $\mat{F}$ are diagonal in a modal Legendre basis with positive
entries, then $\mat{M} \mat{F}[^{-1}] = \mat{F}[^{-1/2}] \mat{M} \mat{F}[^{-1/2}]$ is positive
definite. The modal filters for nodal Gauß- and Lobatto-Legendre as well as modal
Legendre bases described in section \ref{sec:SV-and-filtering} fulfil these
properties.

\begin{lem}
  Augmenting a conservative and stable SBP CPR method for the scalar
  conservation law
  \begin{equation}
  \tag{\ref{eq:scalar-CL}}
    \partial_t u + \partial_x f(u) = 0
  \end{equation}
  with a modal filter $\mat{F}$ applied to the time derivative results in a conservative
  and $\norm{\cdot}_{M F^{-1}}$-stable semidiscretisation if $\mat{M} \mat{F}[^{-1}]$
  induces a scalar product. This condition is fulfilled for nodal bases using
  Gauß- or Lobatto-Legendre points (with lumped mass matrix) and a modal Legendre
  basis.
\end{lem}

This Lemma is connected with the observation of \citex{Allaneau and Jameson}{} \citet{allaneau2011connections},
who presented the reformulation of the energy stable CPR methods given \citex{in}{by}
\citet{vincent2011newclass} as filtered DG methods, where the filter is applied
in the same manner as here. As presented \citex{in}{by} \citet{ranocha2016summation}, these
CPR methods can conserve the discrete norm induced by some positive definite matrix
$\mat{M} + \mat{K}$, corresponding to $\mat{M} \mat{F}[^{-1}]$ in this setting.
However, the discrete norm $\norm{\cdot}_M$ oscillates if the norm
$\norm{\cdot}_{M+K}$ is conserved or even dissipated. Thus, the same behaviour
can be expected here, if the time discretisation is sufficiently accurate.

However, considering the explicit Euler method as time discretisation, the norm
after one time step can be written as
\begin{equation}
  \norm{ \vec{u}_+ }_{M}^2
  = \norm{ \vec{u} + \Delta t \partial_t \vec{u} }_{M}^2
  = \norm{ \vec{u} }_{M}^2
    + 2 \Delta t \left\langle u , \partial_t u \right\rangle_{M}
    + (\Delta t)^2 \norm{ \partial_t \vec{u} }_{M}^2.
\end{equation}
Thus, there is an additional increment of the norm not considered in the
semidiscrete setting of order $(\Delta t)^2$, since the second term is
estimated for semidiscretisations. Applying a filter $\mat{F}$ to the time derivative
yields
\begin{equation}
\begin{aligned}
  \norm{ \vec{u}_+ }_{M F^{-1}}^2
  &=
  \norm{ \vec{u} + \Delta t \mat{F} \partial_t \vec{u} }_{M F^{-1}}^2
  \\&=
  \norm{ \vec{u} }_{M F^{-1}}^2
  + 2 \Delta t \left\langle u , \mat{F} \partial_t u \right\rangle_{M F^{-1}}
  + (\Delta t)^2 \norm{ \mat{F} \partial_t \vec{u} }_{M F^{-1}}^2,
\end{aligned}
\end{equation}
which can be rewritten using $\mat{F}[^T] \mat{M} = \sqrt{\mat{F}[^T]} \mat{M} \sqrt{\mat{F}}$ as
\begin{equation}
  \norm{ \vec{u}_+ }_{M F^{-1}}^2
  = \norm{ \vec{u} }_{M F^{-1}}^2
    + 2 \Delta t \left\langle u , \partial_t u \right\rangle_{M}
    + (\Delta t)^2 \norm{ \mat{F}[^{1/2}] \partial_t \vec{u} }_{M}^2.
\end{equation}
Thus, if the filter $\mat{F}$ reduces the $\mat{M}$-norm, the additional increment
is smaller then in the case of the unfiltered semidiscretisation. By equivalence
of discrete norms, a reduced increase of the norm $\norm{\cdot}_M$ can be
expected.

\subsection{Numerical experiments}
\label{sec:filter-derivative-experiments}

The linear advection equation with discontinuous initial data \eqref{eq:lin-adv-jump}
is solved numerically by an SBP CPR semidiscretisation using $N = 8$ elements
with Gauß-Legendre bases of degree $\leq p = 7$ and an upwind numerical flux
$\fnum(u_-,u_+) = u_-$. The solution on the domain $[0,2]$ with periodic boundary
conditions is advanced in time $t \in [0, 8]$ by an explicit Euler method
using $10^5$ steps.

\begin{figure}[!ht]
  \centering
  \begin{subfigure}[b]{0.49\textwidth}
    \includegraphics[width=\textwidth]{%
      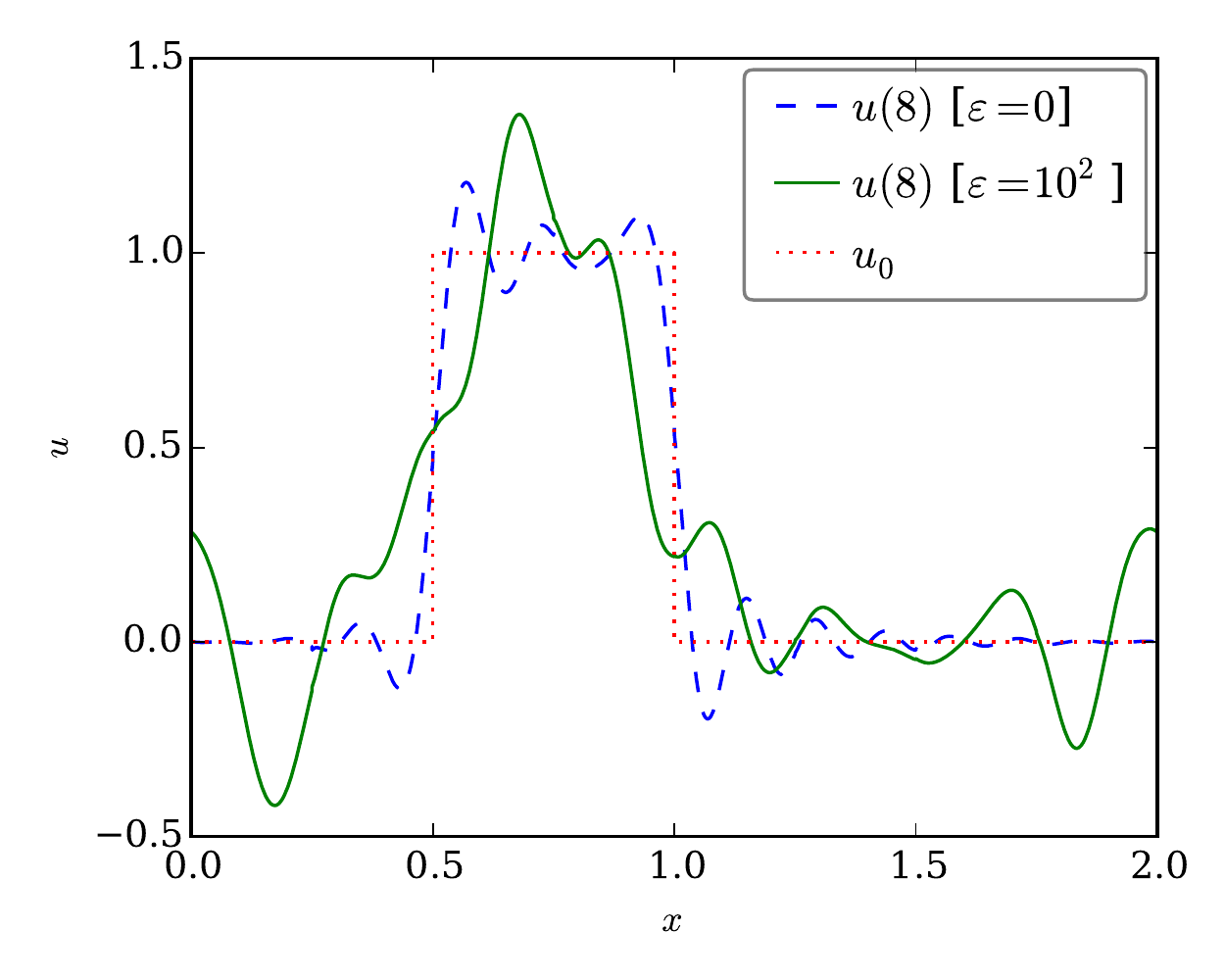}
    \caption{Solution.}
    \label{fig:filtering_time_derivative_solution}
  \end{subfigure}%
  ~
  \begin{subfigure}[b]{0.49\textwidth}
    \includegraphics[width=\textwidth]{%
      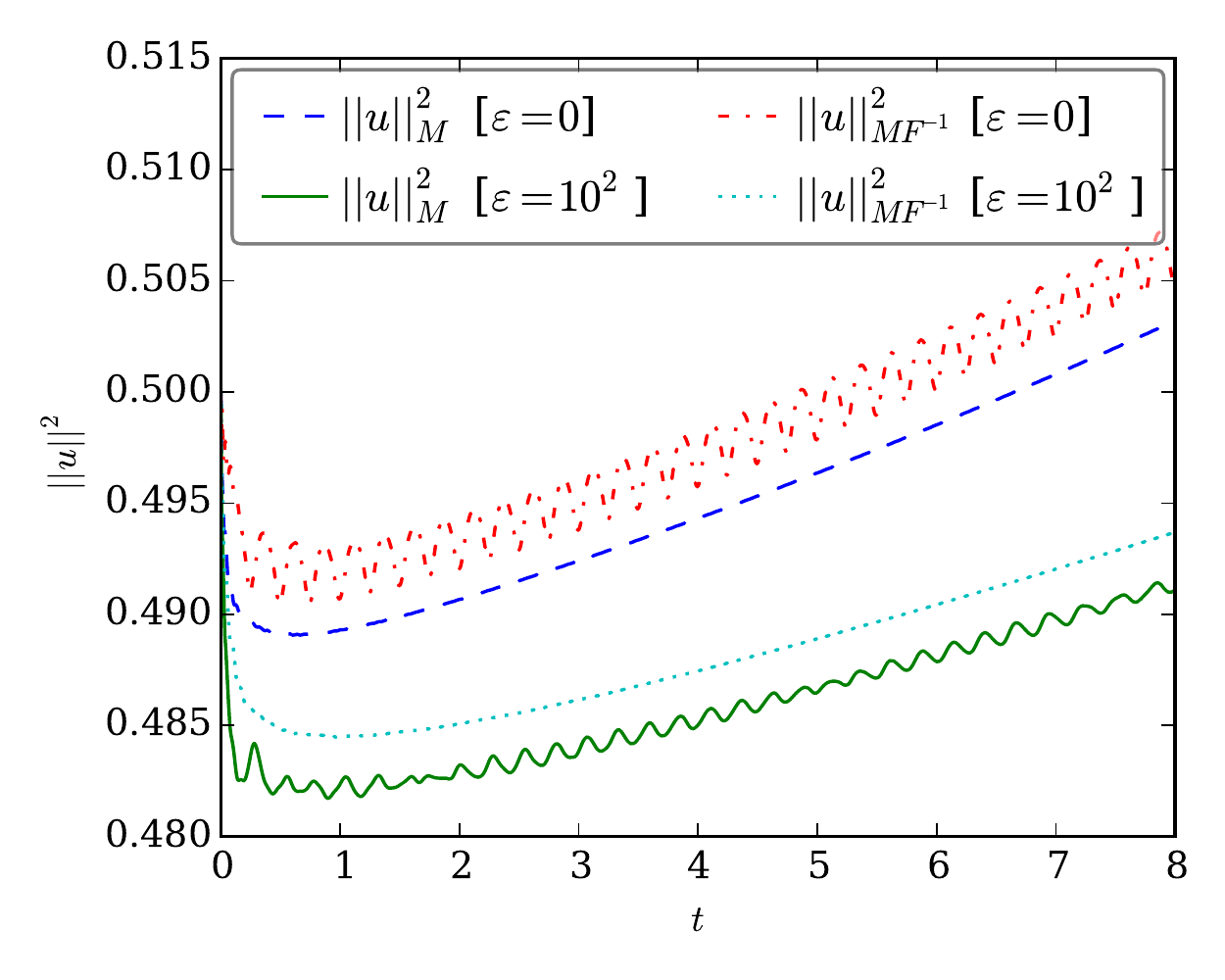}
    \caption{Energy.}
    \label{fig:filtering_time_derivative_energy}
  \end{subfigure}%
  \caption{Numerical results for linear advection using $N = 8$ elements with
           polynomials of degree $\leq p = 7$ with and without modal filtering
           of the time derivative.
           }
  \label{fig:filtering_time_derivative}
\end{figure}

Figure \ref{fig:filtering_time_derivative_solution} shows the solution after
$10^5$ time steps with and without modal filtering of the time derivative,
i.e. for $\epsilon \in \set{0, 100}$. It is obvious that the filter did not reduce
the oscillations. Contrary, these were intensified.

The time development of the norms are plotted in Figure
\ref{fig:filtering_time_derivative_energy}. At first, the norm $\norm{u}_M$ decreases
without filtering. Starting at $t \approx 1$, it increases monotonically.
The corresponding norm $\norm{u}_{M F^{-1}}$ follows the same trend, but is
highly oscillatory. Contrary, the norm $\norm{u}_{M F^{-1}}$ is smoothly
de- and increasing if a modal filter of strength $\epsilon = 100$ is applied to
the time derivative and the corresponding norm $\norm{u}_M$ is oscillatory,
as expected.

  \section{Filtering the solution}

Here, possibility \ref{itm:filter-u} of section \ref{sec:discrete-setting} is
investigated, i.e. the application of a modal filter $\mat{F}$ to the function
$\vec{u}$, used to compute the time derivative.

\subsection{Discrete estimates}

Using again a conservative and stable semidiscretisation of the scalar conservation
law
\begin{equation}
\tag{\ref{eq:scalar-CL}}
  \partial_t u + \partial_x f(u) = 0
\end{equation}
as in section \ref{sec:filter-derivative} results on each element in an ordinary
differential equation
\begin{equation}
  \od{}{t} \vec{u} = g( \vec{u} ).
\end{equation}
Now, a modal filter is applied to the function used to compute the time derivative,
i.e. $\partial_t \vec{u} \mapsto \partial_t \mat{F} \vec{u}$ and the ODE becomes
\begin{equation}
  \od{}{t} \vec{u} = g( \mat{F} \vec{u} ).
\end{equation}

Using an explicit Euler method as time discretisation, the value after one
time step of size $\Delta t$ can be written as
\begin{equation}
  \vec{u}_+ 
  = \vec{u} + \Delta t \, \partial_t \mat{F} \vec{u}.
\end{equation}
If the filter $\mat{F}$ is invertible, this can be rewritten as
\begin{equation}
\label{eq:rewritten-filter-solution}
  \vec{u}_+ 
  = \mat{F}[^{-1}] \left( \mat{F} \vec{u}
    + \Delta t \, \mat{F} \partial_t \mat{F} \vec{u} \right).
\end{equation}
The term in brackets corresponds to a discrete combination of the approaches
\ref{itm:operator-splitting} and \ref{itm:filter-derivative} mentioned in section
\ref{sec:discrete-setting}: At first, the filter $\mat{F}$ is applied in a split operator
fashion to the solution $\vec{u}$ and afterwards a time step using a scheme
with filter $\mat{F}$ applied to the time derivative is used.

If the assumptions of the previous sections are complied with, the second step
corresponds to a stable semidiscretisation with respect to the norm induced by
$\mat{M} \mat{F}[^{-1}]$. The additional application of a filter prior to this
step could control the additional terms similarly to the way described in section
\ref{sec:operator-splitting}.

However, after a full time step of this split operator formulation with filtered
time derivative, the inverse filter $\mat{F}[^{-1}]$ is applied, destroying this
stability estimates. Therefore, we do not expect this scheme to be superior
compared to the other possibilities mentioned in section \ref{sec:discrete-setting}.

\subsection{Numerical experiments}

The same setup as in the previous section \ref{sec:filter-derivative-experiments}
is used to compute numerical solutions of the linear advection equation with
discontinuous initial data \eqref{eq:lin-adv-jump}.

\begin{figure}[!ht]
  \centering
  \begin{subfigure}[b]{0.49\textwidth}
    \includegraphics[width=\textwidth]{%
      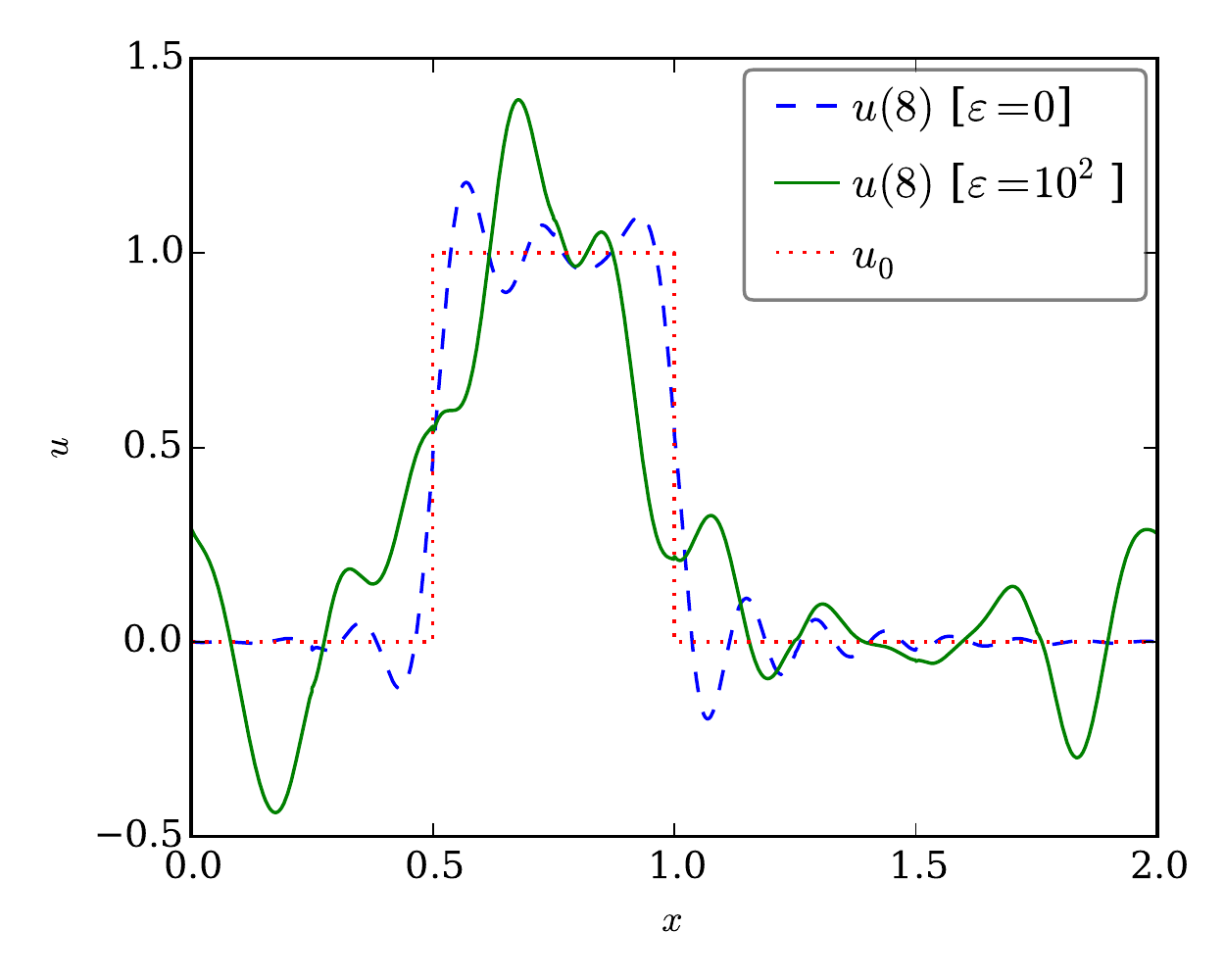}
    \caption{Solution.}
    \label{fig:filtering_solution_solution}
  \end{subfigure}%
  ~
  \begin{subfigure}[b]{0.49\textwidth}
    \includegraphics[width=\textwidth]{%
      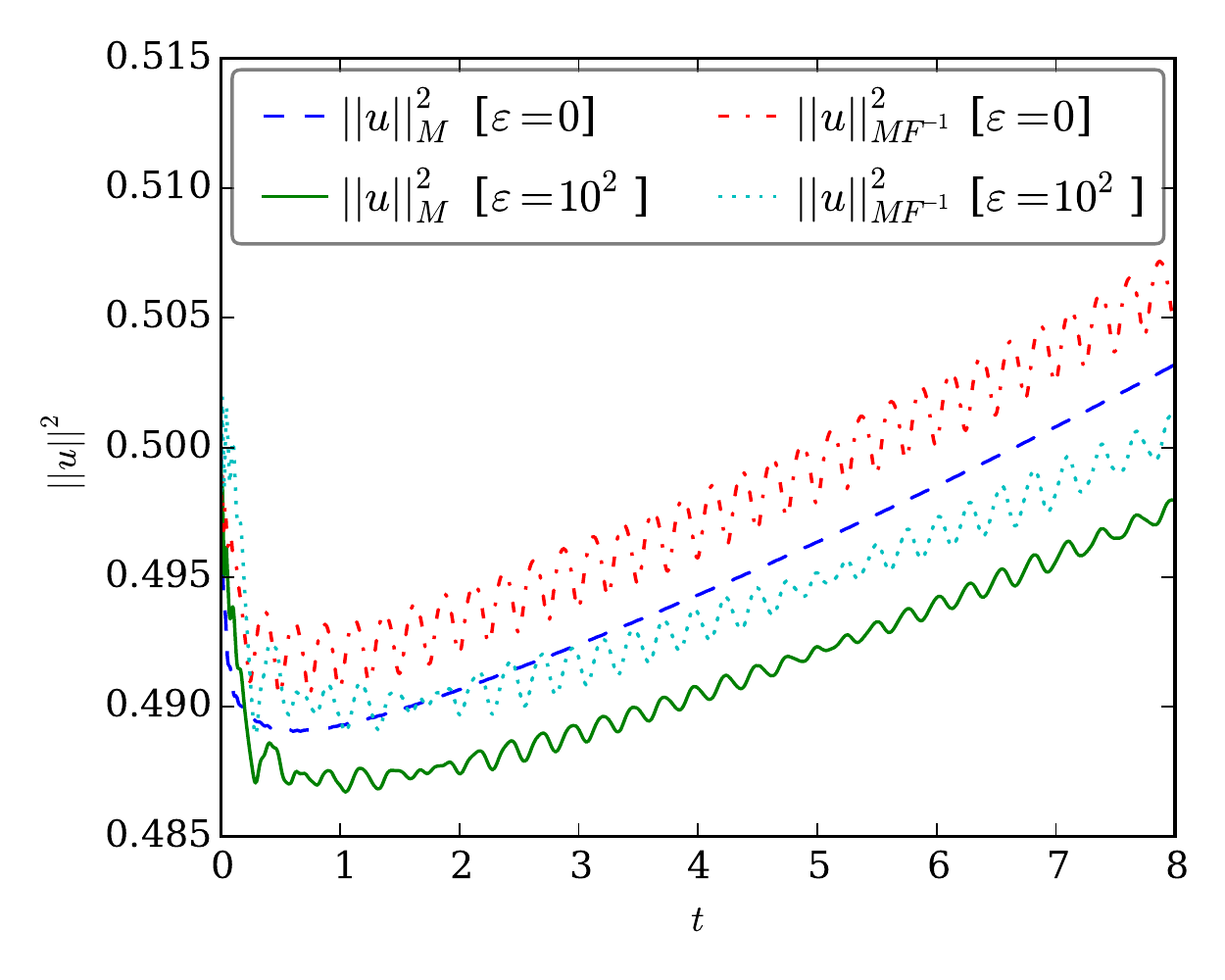}
    \caption{Energy.}
    \label{fig:filtering_solution_energy}
  \end{subfigure}%
  \caption{Numerical results for linear advection using $N = 8$ elements with
           polynomials of degree $\leq p = 7$ with and without modal filtering
           of the solution.
           }
  \label{fig:filtering_solution}
\end{figure}

The solution in Figure \ref{fig:filtering_solution_solution} is very similar to
the one in Figure \ref{fig:filtering_time_derivative_solution} of the previous
section, as could be guessed from equation \eqref{eq:rewritten-filter-solution},
if the first application of the filter $\mat{F}$ and the inverse filter
$\mat{F}[^{-1}]$ are ignored. Especially, the application of modal filtering in
this setting does not decrease oscillations or increase the quality of the
solution.

However, the influence of the additional terms in equation
\eqref{eq:rewritten-filter-solution} become visible in the energy of the numerical
solution in Figure \ref{fig:filtering_solution_energy}. The norms without
application of modal filtering are the same as in Figure
\ref{fig:filtering_time_derivative_energy} of the previous section, whereas
the norms $\norm{u}_M$ and $\norm{u}_{M F^{-1}}$ of the filtered solution are
both oscillating.

  \section{Conclusions and further research}

In this work, modal filtering has been applied in the general framework of
CPR methods using SBP operators. It is shown to be equivalent to the application
of artificial dissipation / spectral viscosity up to first order in time
and circumvents the time step restrictions of these schemes described in the
first part of this series \citex{}{by} \citet{ranocha2016enhancing}.

Additionally, a new adaptive method to chose the filter strength automatically
has been proposed. Compensating error terms of order $(\Delta t)^2$, the
fully discrete schemes using an explicit Euler method become provably stable.
However, this adaptive filtering does not remove all oscillations of the numerical
solutions, especially
in the nonlinear case developing discontinuities.

Moreover, two other possibilities of the application of modal filtering are
investigated. Theoretical considerations about their inferior stability
properties are accompanied by numerical experiments verifying these.

However, the analysis has been limited to the explicit Euler method as time
discretisation. Although these results carry over to strong-stability preserving
methods composed of explicit Euler steps, specialised investigations will be
conducted.

Of course, an extension of these results to other hyperbolic conservation laws
will be interesting.

  \printbibliography

\end{document}